\documentclass[12pt,twoside]{article}
\usepackage[english]{babel}
\usepackage[latin1]{inputenc}
\usepackage{amsmath}
\usepackage{amssymb,amsfonts}
\usepackage{graphicx}                   

\newcommand{\bR}{\mathbf{R}}

\newcommand{\bb}{\mathbf{b}}

\newcommand{\bt}{\mathbf{t}}

\newcommand{\HYP}{\mathbb{H}^3}
\newcommand{\HYN}{\mathbb{H}^n}

\begin{document}
\pagestyle{myheadings}
\markboth{\centerline{Jen\H o Szirmai}}
{Hyperball packings related to octahedron and cube tilings $\dots$ }
\title
{Hyperball packings related to octahedron and cube tilings in hyperbolic space}

\author{\normalsize{Jen\H o Szirmai} \\
\normalsize Budapest University of Technology and \\
\normalsize Economics Institute of Mathematics, \\
\normalsize Department of Geometry \\
\date{\normalsize{\today}}}

\maketitle


\begin{abstract}
In this paper we study congruent and non-congruent hyperball (hypersphere) packings of the truncated regular octahedron and cube tilings.
These are derived from the Coxeter simplex tilings $\{p,3,4\}$ $(7\le p \in \mathbb{N})$ and $\{p,4,3\}$ $(5\le p \in \mathbb{N})$
in $3$-dimensional hyperbolic space $\HYP$.
We determine the densest hyperball packing arrangement and its density
with congruent and non-congruent hyperballs related to the above tilings in $\HYP$.

We prove that the locally densest congruent or non-congruent hyperball configuration belongs to the regular truncated cube with density
$\approx 0.86145$. This is larger than the B\"or\"oczky-Florian density upper bound for balls and horoballs.
Our locally optimal non-congruent hyperball packing configuration cannot be extended to the entire 
hyperbolic space $\mathbb{H}^3$, but we determine the extendable densest non-congruent hyperball packing arrangement related to a regular cube tiling 
with density $\approx 0.84931$.

\end{abstract}

\newtheorem{theorem}{Theorem}[section]
\newtheorem{corollary}[theorem]{Corollary}
\newtheorem{conjecture}{Conjecture}[section]
\newtheorem{lemma}[theorem]{Lemma}
\newtheorem{exmple}[theorem]{Example}
\newtheorem{defn}[theorem]{Definition}
\newtheorem{rmrk}[theorem]{Remark}
\newenvironment{definition}{\begin{defn}\normalfont}{\end{defn}}
\newenvironment{remark}{\begin{rmrk}\normalfont}{\end{rmrk}}
\newenvironment{example}{\begin{exmple}\normalfont}{\end{exmple}}
\newenvironment{acknowledgement}{Acknowledgement}


\section{Introduction}
In $n$-dimensional hyperbolic space $\mathbb{H}^n$ $(n\ge2)$ there are $3$ kinds
of ''balls (spheres)": the classical balls (spheres), horoballs (horospheres) and hyperballs (hyperspheres).

In this paper we consider the hyperballs and their packings in $3$-dimensional hyperbolic space $\HYP$. However, first we survey the previous results related to this topic.

In the hyperbolic plane $\mathbb{H}^2$ the universal upper bound of the hypercycle packing density is $\frac{3}{\pi}$,
proved by I.~Vermes in \cite{V79} and the universal lower bound of the hypercycle covering density is $\frac{\sqrt{12}}{\pi}$
determined by I.~Vermes in \cite{V81}. 

In \cite{Sz06-1} and \cite{Sz06-2} we analysed the regular prism tilings (simple truncated Coxeter orthoscheme tilings) and the corresponding optimal hyperball packings in
$\mathbb{H}^n$ $(n=3,4)$ and we extended the method developed in the former paper \cite{Sz13-3} to 
5-dimensional hyperbolic space.
In paper \cite{Sz13-4} we studied the $n$-dimensional hyperbolic regular prism honeycombs
and the corresponding coverings by congruent hyperballs and we determined their least dense covering densities.
Furthermore, we formulated conjectures for the candidates of the least dense hyperball
covering by congruent hyperballs in the 3- and 5-dimensional hyperbolic space ($n \in \mathbb{N},3 \le n \le 5)$.

In \cite{Sz17-1} we discussed congruent and non-congruent hyperball packings of the truncated regular tetrahedron tilings.
These are derived from the Coxeter simplex tilings $\{p,3,3\}$ $(7\le p \in \mathbb{N})$ and $\{5,3,3,3,3\}$
in $3$- and $5$-dimensional hyperbolic space.
We determined the densest hyperball packing arrangement and its density
with congruent hyperballs in $\mathbb{H}^5$ and determined the smallest density upper bounds of 
non-congruent hyperball packings generated by the above tilings in $\HYN,~ (n=3,5)$.

In \cite{Sz17} we deal with the packings derived by horo- and hyperballs (briefly hyp-hor packings) in $n$-dimensional hyperbolic spaces $\HYN$
($n=2,3$) which form a new class of the classical packing problems.
We constructed in the $2-$ and $3-$dimensional hyperbolic spaces hyp-hor packings that
are generated by complete Coxeter tilings of degree $1$ 
and we determined their densest packing configurations and their densities.
We proved using also numerical approximation methods that in the hyperbolic plane ($n=2$) the density of the above hyp-hor packings arbitrarily approximate
the universal upper bound of the hypercycle or horocycle packing density $\frac{3}{\pi}$ and
in $\HYP$ the optimal configuration belongs to the $\{7,3,6\}$ Coxeter tiling with density $\approx 0.83267$.
Furthermore, we analyzed the hyp-hor packings in
truncated orthosche\-mes $\{p,3,6\}$ $(6< p < 7, ~ p\in \mathbb{R})$ whose
density function is attained its maximum for a parameter which lies in the interval $[6.05,6.06]$
and the densities for parameters lying in this interval are larger that $\approx 0.85397$. That means that these
locally optimal hyp-hor configurations provide larger densities that the B\"or\"oczky-Florian density upper bound
$(\approx 0.85328)$ for ball and horoball packings but these hyp-hor packing configurations cannot be extended to the entirety of hyperbolic space $\mathbb{H}^3$. 

In \cite{Sz14} we studied a large class of hyperball packings in $\HYP$
that can be derived from truncated tetrahedron tilings (see e.g. \cite{S14}, \cite{S17}).
We proved that if the truncated tetrahedron is regular, then the density
of the densest packing is $\approx 0.86338$. This is larger than the B\"or\"oczky-Florian density upper bound
but our locally optimal hyperball packing configuration cannot be extended to the entirety of
$\mathbb{H}^3$. However, we described a hyperball packing construction, 
by the regular truncated tetrahedron tiling under the extended Coxeter group $\{3, 3, 7\}$ with maximal density $\approx 0.82251$.

Recently, (to the best of author's knowledge) the candidates for the densest hyperball
(hypersphere) packings in the $3,4$ and $5$-dimensional hyperbolic space $\mathbb{H}^n$ are derived by the regular prism
tilings which have been in papers \cite{Sz06-1}, \cite{Sz06-2} and \cite{Sz13-3}.

{In \cite{Sz17-2} we considered hyperball packings in 
$3$-dimensional hyperbolic space. We developed a decomposition algorithm that for each saturated hyperball packing provides a decomposition of $\HYP$ 
into truncated tetrahedra. Therefore, in order to get a density upper bound for hyperball packings, it is sufficient to determine
the density upper bound of hyperball packings in truncated simplices.}

{\it Now, we consider packings related to truncated regular octahedron and cube tilings that are derived from the Coxeter simplex tilings 
$\{p,3,4\}$ $(7\le p \in \mathbb{N})$ and $\{p,4,3\}$ $(5\le p \in \mathbb{N})$
in $3$-dimensional hyperbolic space $\HYP$.
We determine the densest hyperball packing arrangement and its density
with congruent and non-congruent hyperballs related to the above tilings. 
Moreover, we prove that the locally densest congruent or non-congruent hyperball configuration belongs to the regular truncated cube with density
$\approx 0.86145$. This is larger than the B\"or\"oczky-Florian density upper bound for balls and horoballs.
Our locally optimal non-congruent hyperball packing configuration cannot be extended to the entire 
hyperbolic space $\mathbb{H}^3$, but we describe a non-congruent hyperball packing construction, 
by the regular cube tiling under the extended Coxeter group $\{4, 3, 7\}$ 
with maximal density $\approx 0.84931$.
The main results are summarized in Theorems 3.2, 3.6 and in Corollary 3.3 related to the octahedron tilings, moreover, related to the cube tilings in 
Theorems 3.8, 3.12-15 and in Corollary 3.9.}
\section{Basic notions}
We use for $\mathbb{H}^3$ (and analogously for $\HYN$, $n\ge3$) the projective model
in the Lorentz space $\mathbb{E}^{1,3}$
that denotes the real vector space $\mathbf{V}^{4}$ equipped with the bilinear
form of signature $(1,3)$,
$
\langle \mathbf{x},~\mathbf{y} \rangle = -x^0y^0+x^1y^1+x^2y^2+ x^3 y^3,
$
where the non-zero vectors
$
\mathbf{x}=(x^0,x^1,x^2,x^3)\in\mathbf{V}^{4} \ \  \text{and} \ \ \mathbf{y}=(y^0,y^1,y^2,y^3)\in\mathbf{V}^{4},
$
are determined up to real factors, for representing points 
of $\mathcal{P}^3(\mathbb{R})$. Then $\mathbb{H}^3$ can be interpreted
as the interior of the conical quadric
$
Q=\{(\mathbf{x})\in\mathcal{P}^3 | \langle  \mathbf{x},~\mathbf{x} \rangle =0 \}=:\partial \mathbb{H}^3
$
in the real projective space $\mathcal{P}^3(\mathbf{V}^{4},
\mbox{\boldmath$V$}\!_{4})$ (here $\mbox{\boldmath$V$}\!_{4}$ is the dual space of $\mathbf{V}^{4}$). 
Namely, for an interior point $\mathbf{y}$ there holds $\langle  \mathbf{y},~\mathbf{y} \rangle <0$.
(Restricting this model to the hyperplane $x^0=1$ we obtain the usual collinear, i.e., Cayley-Klein model.)

Points of the boundary $\partial \mathbb{H}^3 $ in $\mathcal{P}^3$
are called points at infinity, or at the absolute of $\mathbb{H}^3 $. Points lying outside $\partial \mathbb{H}^3 $
are said to be outer points of $\mathbb{H}^3 $ relative to $Q$. Let $(\mathbf{x}) \in \mathcal{P}^3$, a point
$(\mathbf{y}) \in \mathcal{P}^3$ is said to be conjugate to $(\mathbf{x})$ relative to $Q$ if
$\langle \mathbf{x},~\mathbf{y} \rangle =0$ holds. The set of all points which are conjugate to $(\mathbf{x})$
form a projective (polar) hyperplane
$
pol(\mathbf{x}):=\{(\mathbf{y})\in\mathcal{P}^3 | \langle \mathbf{x},~\mathbf{y} \rangle =0 \}.
$
Thus, the quadric $Q$ induces a bijection
(linear polarity $\mathbf{V}^{4} \rightarrow
\mbox{\boldmath$V$}\!_{4})$
from the points of $\mathcal{P}^3$ onto their polar hyperplanes.

Point $X (\bold{x})$ and hyperplane $\alpha (\mbox{\boldmath$a$})=\{(x^0,x^1,x^2,x^3)|\sum_{i=0}^3 x^i a_i=0 \}$
are incident if $\bold{x}\mbox{\boldmath$a$}=0$ ($\bold{x} \in \bold{V}^{4} \setminus \{\mathbf{0}\}, \ \mbox{\boldmath$a$}
\in \mbox{\boldmath$V$}_{4}
\setminus \{\mbox{\boldmath$0$}\}$).

The hypersphere (or equidistant surface) is a quadratic surface at a constant distance
from a plane (base plane) in both halfspaces. The infinite body bounded by the hypersphere, containing the base plane, is called {\it hyperball}.

The {\it half hyperball } (i.e., the part of the hyperball lying on one side of its base plane) with distance $h$ to a base plane $\beta$ is denoted by $\mathcal{H}^h_+$.
The volume of the intersection of $\mathcal{H}^h_+(\mathcal{A})$ and the right prism with
base a $2$-polygon $\mathcal{A} \subset \beta$ can be determined by the classical formula 
(2.1) of J.~Bolyai \cite{B91}.
\begin{equation}
\mathrm{Vol}(\mathcal{H}^h_+(\mathcal{A}))=\frac{1}{4}\mathrm{Area}(\mathcal{A})\left[k \sinh \frac{2h}{k}+
2 h \right], \tag{2.1}
\end{equation}
The constant $k =\sqrt{\frac{-1}{K}}$ is the natural length unit in
$\mathbb{H}^3$, where $K$ denotes the constant negative sectional curvature. In the following we may assume that $k=1$.
\subsection{Complete orthoschemes}
\begin{definition}
An orthoscheme $\mathcal{S}$ in $\mathbb{H}^n$ $(2\le n \in \mathbb{N})$ is a simplex bounded by $n+1$ hyperplanes $H^0,\dots,H^n$
such that
(see \cite{K89})
$
H^i \bot H^j, \  \text{for} \ j\ne i-1,i,i+1.
$
\end{definition}

{\it The orthoschemes of degree} $m$ in $\mathbb{H}^n$ are bounded by $n+m+1$ hyperplanes
$H^0,H^1,\dots,H^{n+m}$ such that $H^i \perp H^j$ for $j \ne i-1,~i,~i+1$, where, for $m=2$,
indices are taken modulo $n+3$. For a usual (classical) orthoscheme we denote the $(n+1)$-hyperface opposite to the vertex $A_i$
by $H^i$ $(0 \le i \le n)$. An orthoscheme $\mathcal{S}$ has $n$ dihedral angles which
are not right angles. Let $\alpha^{ij}$ denote the dihedral angle of $\mathcal{S}$
between the faces $H^i$ and $H^j$. Then we have
$
\alpha^{ij}=\frac{\pi}{2}, \ \ \text{if} \ \ 0 \le i < j -1 \le n.
$
The $n$ remaining dihedral angles $\alpha^{i,i+1}, \ (0 \le i \le n-1)$ are called the
essential angles of $\mathcal{S}$.
Geometrically, complete orthoschemes of degree $d$ can be described as follows:
\begin{enumerate}
\item
For $m=0$, they coincide with the class of classical orthoschemes introduced by
{{Schl\"afli}} (see Definitions 2.1).
The initial and final vertices, $A_0$ and $A_n$ of the orthogonal edge-path
$A_iA_{i+1},~ i=0,\dots,n-1$, are called principal vertices of the orthoscheme.
\item
A complete orthoscheme of degree $m=1$ can be interpreted as an
orthoscheme with one outer principal vertex, say $A_n$, which is truncated by
its polar plane $pol(A_n)$ (see Fig.~1 and 3). In this case the orthoscheme is called simply truncated with
outer vertex $A_n$.
\item
A complete orthoscheme of degree $m=2$ can be interpreted as an
orthoscheme with two outer principal vertices, $A_0,~A_n$, which is truncated by
its polar hyperplanes $pol(A_0)$ and $pol(A_n)$. In this case the orthoscheme is called doubly
truncated. We distinguish two different types of orthoschemes but 
will not enter into the details (see \cite{K89}).
\end{enumerate}

The ordered set $\{k_1,\dots,k_{n-1},k_n\}$ is said to be the
Coxeter-Schl\"afli symbol of the simplex tiling $\mathcal{P}$ generated by $\mathcal{S}$.
To every scheme there is a corresponding
symmetric matrix $(c^{ij})$ of size $(n+1)\times(n+1)$ where $c^{ii}=1$ and, for $i \ne j\in \{0,1,2,\dots,n \}$,
$c^{ij}$ equals $-\cos{\frac{\pi}{k_{ij}}}$ with all angles between the facets $i$,$j$ of $\mathcal{S}$.

For example, $(c^{ij})$ below is the so called Coxeter-Schl\"afli matrix of the orthoscheme $S$ in
3-dimensional hyperbolic space $\mathbb{H}^3$ with
parameters (nodes) $k_1=p,k_2=q,k_3=r$ :
\[
(c^{ij}):=\begin{pmatrix}
1& -\cos{\frac{\pi}{p}}& 0 & 0 \\
-\cos{\frac{\pi}{p}} & 1 & -\cos{\frac{\pi}{q}}& 0 \\
0 & -\cos{\frac{\pi}{q}} & 1 & -\cos{\frac{\pi}{r}} \\
0 & 0 & -\cos{\frac{\pi}{r}} & 1 \\
\end{pmatrix}. \tag{2.2}
\]
In general, the complete Coxeter orthoschemes were classified by {{Im Hof}} in
\cite{IH85} and \cite{IH90} by generalizing the method of {{Coxeter}} and {{B\"ohm}}, who
showed that they exist only for dimensions $\leq 9$. From this classification it follows, that the complete
orthoschemes of degree $m=1$ exist up to 5 dimensions.

In this paper we consider some tilings generated by orthoschemes of degree $1$ where the initial vertex $A_n$ is an outer point regarding the quadric $Q$. 
These orthoschemes and the corresponding Coxeter tilings exist in the $2$-, $3-$, $4-$ and
$5-$dimensional hyperbolic spaces and
are characterized by their Coxeter-Schl\"afli symbols and graphs.

In $n$-dimensional hyperbolic space $\mathbb{H}^n$ $(n \ge 2)$
it can be seen that if $\mathcal{S}=A_0A_1A_2 \dots A_n$ $P_0P_1P_2 \dots P_n$ is a complete
orthoscheme  with degree $m=1$ (a simply frustum orthoscheme) where $A_n$ is a outer vertex of
$\mathbb{H}^n$ then the points $P_0,P_1,P_2,\dots,P_{n-1}$ lie on the polar hyperplane $\pi$ of $A_n$ (see Fig.~1 in $\HYP$).
\begin{figure}[ht]
\centering
\includegraphics[width=6.5cm]{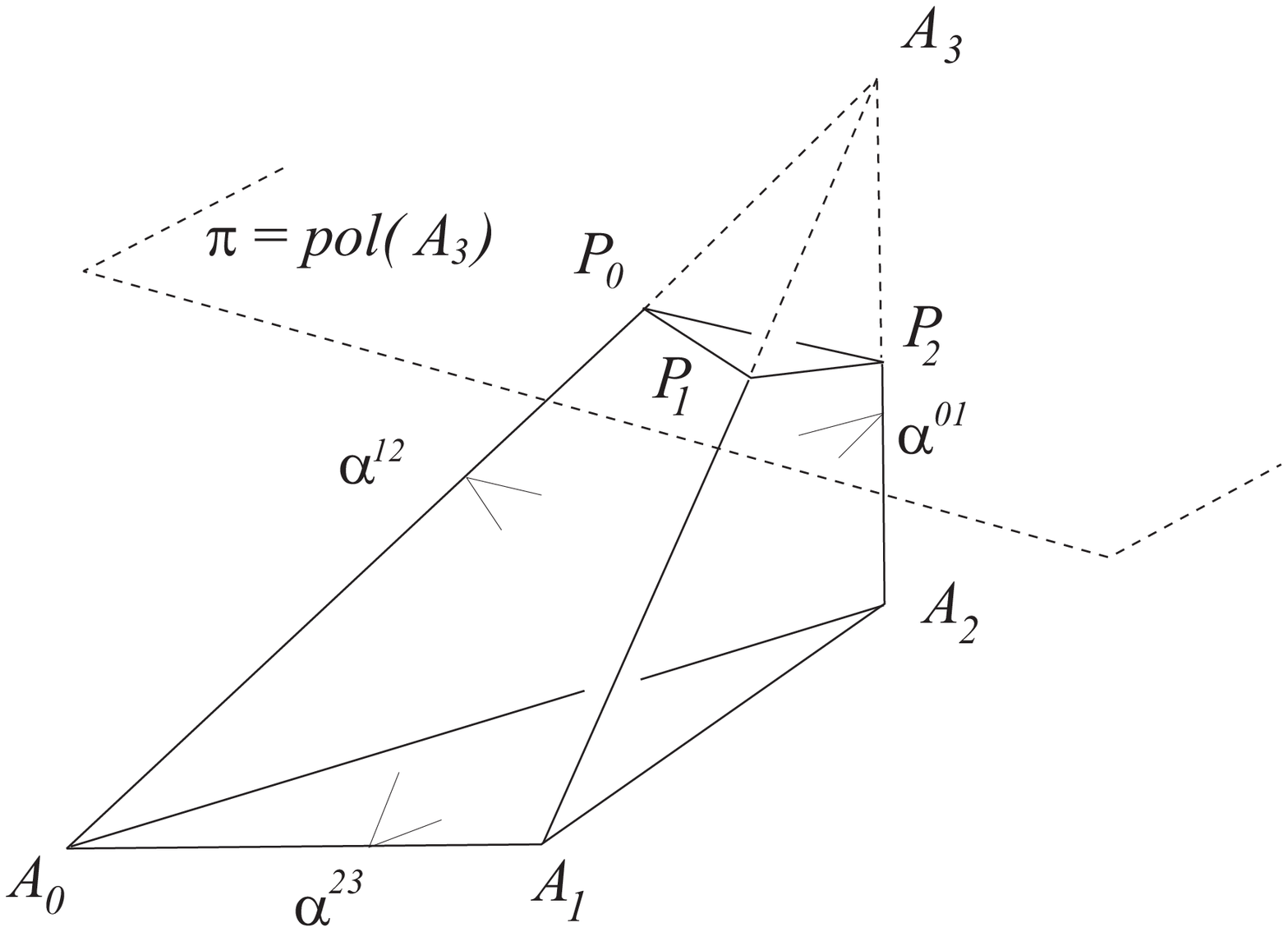} \includegraphics[width=6.5cm]{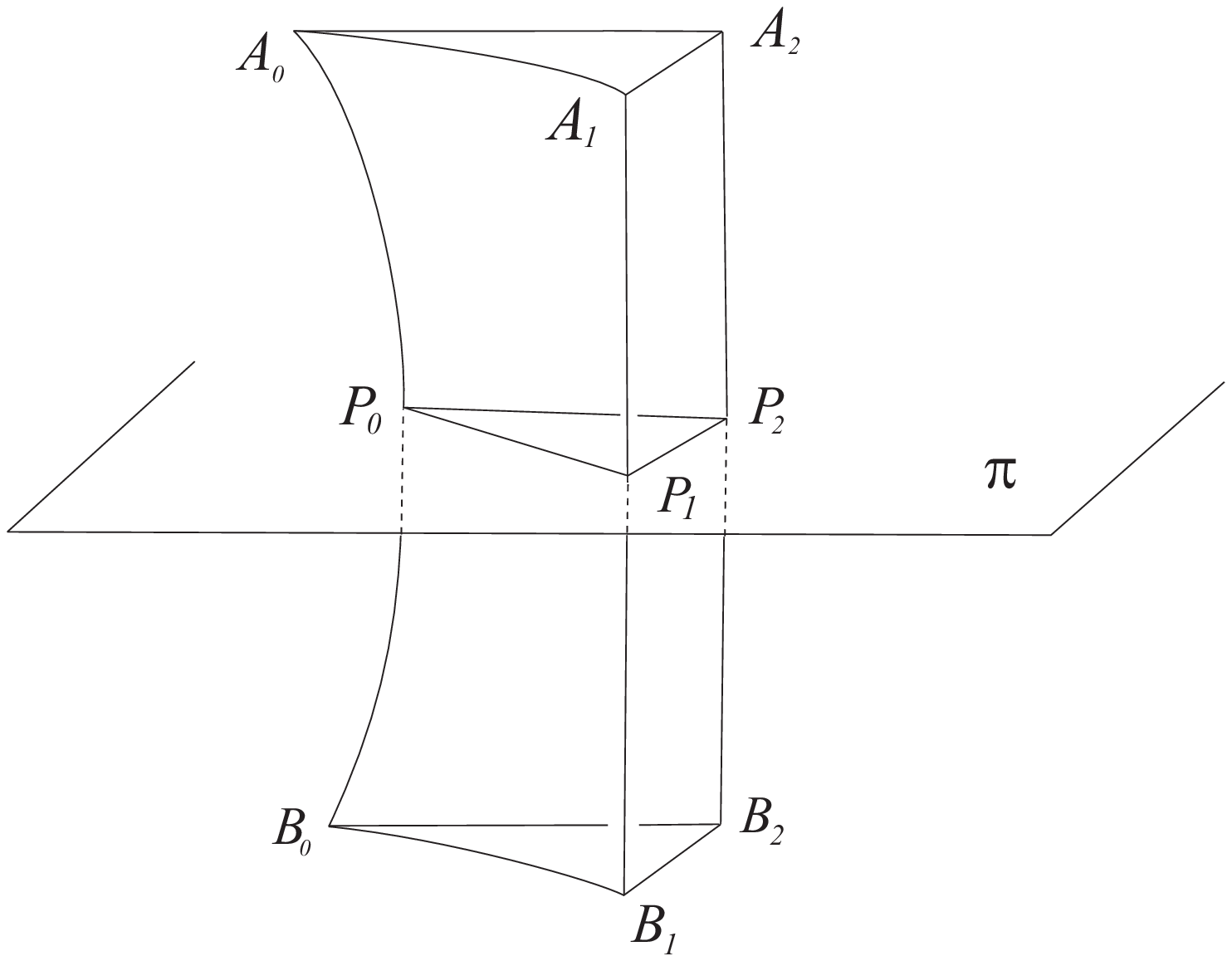}

a. \hspace{5cm} b.
\caption{a.~A $3$-dimensional complete orthoscheme of degree $m=1$ (simple frustum orthoscheme)
with outer vertex $A_3$. This orthoscheme is truncated by its polar plane $\pi=pol(A_3)$.~b.~Two congruent
adjacent simple frustum orthoschemes.}
\label{}
\end{figure}
The images of $\mathcal{O}$ under reflections on its side facets generate a tiling in $\HYN$.
{Our polyhedron $A_0A_1A_2P_0P_1P_2$ is a simple frustum orthoscheme with
outer vertex $A_3$ (see Fig.~1) whose volume can be calculated by the following theorem of R.~Kellerhals
\cite{K89}:}
\begin{theorem} The volume of a three-dimensional hyperbolic
complete ortho\-scheme (except Lambert cube cases) $\mathcal{S}$
is expressed with the essential angles $\alpha_{01},\alpha_{12},\alpha_{23}, \ (0 \le \alpha_{ij} \le \frac{\pi}{2})$
(Fig.~1) in the following form:
\begin{align}
&Vol_3(\mathcal{O})=\frac{1}{4} \{ \mathcal{L}(\alpha_{01}+\theta)-
\mathcal{L}(\alpha_{01}-\theta)+\mathcal{L}(\frac{\pi}{2}+\alpha_{12}-\theta)+ \notag \\
&+\mathcal{L}(\frac{\pi}{2}-\alpha_{12}-\theta)+\mathcal{L}(\alpha_{23}+\theta)-
\mathcal{L}(\alpha_{23}-\theta)+2\mathcal{L}(\frac{\pi}{2}-\theta) \}, \tag{2.3}
\end{align}
where $\theta \in [0,\frac{\pi}{2})$ is defined by the following formula:
$$
\tan(\theta)=\frac{\sqrt{ \cos^2{\alpha_{12}}-\sin^2{\alpha_{01}} \sin^2{\alpha_{23}
}}} {\cos{\alpha_{01}}\cos{\alpha_{23}}}
$$
and where $\mathcal{L}(x):=-\int\limits_0^x \log \vert {2\sin{t}} \vert dt$ \ denotes the
Lobachevsky function.
\end{theorem}
For our prism tilings $\mathcal{T}_{pqr}$ we have:~
$\alpha_{01}=\frac{\pi}{p}, \ \ \alpha_{12}=\frac{\pi}{q}, \ \
\alpha_{23}=\frac{\pi}{r}$ .
\section{On hyperball packings related to truncated octahedron and cube tilings}
{Similarly to the truncted tetrahedral cases (see \cite{Sz06-1}, \cite{Sz06-2}, \cite{Sz13-4}, \cite{Sz14}, \cite{Sz17-1}, \cite{Sz17-2}) it is intersting to examine and 
to construct locally optimal {\it congruent and non-congruent} hyperball packings and coverings
relating to suitable truncated polyhedron tilings in $3$- and higher dimensions as well.

{\it In this paper we consider the 3-dimensional regular, truncated octahedron and cube tilings that are derived 
from the Coxeter simplex tilings $\{p_1,4,3\}$ $(\mathbb{N} \ni p_1 \ge 5)$ and $\{p_2,3,4\}$ $(\mathbb{N} \ni p_2 \ge 7)$.}

\subsection{Hyperball packings with congruent hyperballs related to regular truncated octahedron tiling $\{3,4,p\}$}
We consider a regular truncated octahedron tiling $\mathcal{T}(\mathcal{O}^r(p))$ with Schl\"afli symbol $\{3,4,p\}$, $(5 \le p\in\mathbb{N})$.
These tilings are derived by duality from the Coxeter tilings $\{p,4,3\}$ whose fundamental 
domains are simply truncated orthoschems (e.g. $P_1P_2P_3Q_1Q_2Q_3$ in Fig.~2.b). 

Let a truncated octahedra $\mathcal{O}^r(p)$ $\subset \HYP$ be a tile from the above tiling 
that is illustrated in Fig.~2.a.b. This truncated octahedra can be derived also by truncation from a regular Euclidean octahedron centred at the origin with vertices $B_i$
($i\in\{1,\dots,6\}$. The truncating planes $\beta_i$ are the polar planes of outer vertices $B_i$ that can be the ultraparallel base planes 
of hyperballs with height $s$ $\mathcal{H}^{s}_i$ $(i=1,\dots,6)$. The non-orthogonal dihedral angles of $\mathcal{O}^r(p)$ are equal to $\frac{2\pi}{p}$.
The distances between two base planes are equal $d(\beta_i,\beta_j)=:e_{ij}=2h(p)$ ($i < j$, $i,j \in \{1,\dots 6\})$ ($d$ is the hyperbolic distance function)
therefore the height of a hyperball is at most $h(p)$ (see Fig.~2.a.b). It is clear, that in the optimal case the heights of the hyperballs are $h(p)$ 
i.e. the congruent hyperballs touch each other.

We consider a saturated congruent hyperball packing $\mathcal{B}^{h(p)}$ of hyperballs $\mathcal{H}^{h(p)}$ related to the above regular, truncated octahedron tilings 
$\{3,4,p\}$ $(\mathbb{N} \ni p \ge 5)$.
\begin{figure}[ht]
\centering
\includegraphics[width=6.5cm]{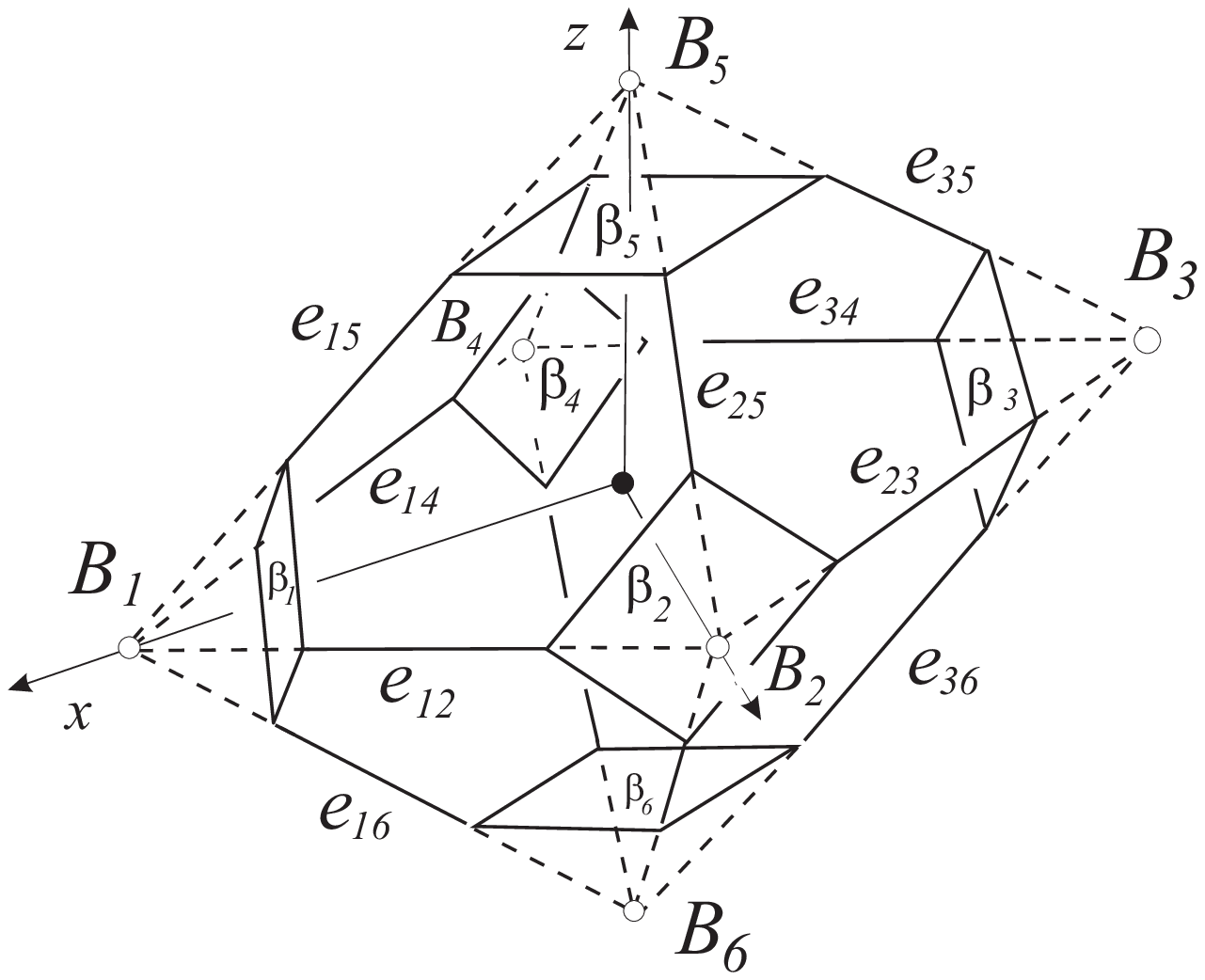} \includegraphics[width=6.5cm]{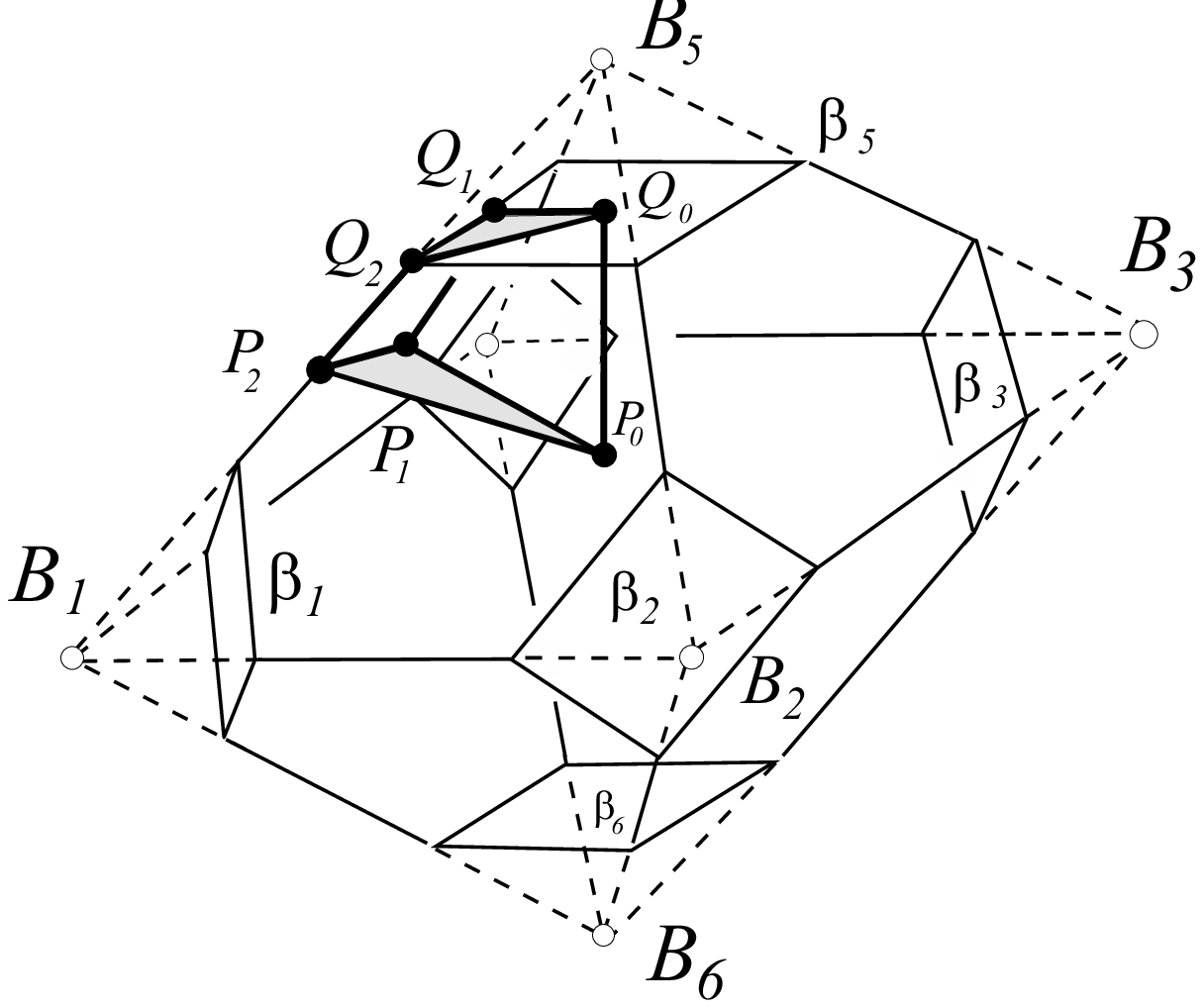}

a. \hspace{6cm} b.
\caption{Regular truncated octahedron. ~b.~Regular truncated octahedron with complete orthoscheme of degree $m=1$ (simple frustum orthoscheme)
with outer vertex $B_5$. This orthoscheme is truncated by its polar plane $\beta_5=pol(B_5)$.}
\label{}
\end{figure}
The volume of the truncated octahedron $\mathcal{O}^r(p)$ is denoted by $Vol(\mathcal{O}^r(p))$ and
we introduce the locally density function $\delta(\mathcal{O}^r(h(p)))$ related to $\mathcal{O}^r(p)$:
\begin{definition}
\begin{equation}
\delta(\mathcal{O}^r(h(p))):=\frac{6\cdot Vol(\mathcal{H}^{h(p)} \cap \mathcal{O}^r(p))}{Vol({\mathcal{O}^r(p)})}. \tag{3.1}
\end{equation}
\end{definition}
If the parameter $p$ is given then the common length of the common perpendiculars $2h(p)=e_{ij}$ $(i < j$, $i,j \in \{1,\dots,6\})$
can be determined by the machinery of the projective geometry.

The points $P_2[{\mathbf{p}}_2]$ and $Q_2[{\mathbf{q}}_2]$ are proper points of hyperbolic $3$-space and
$Q_2$ lies on the polar hyperplane $pol(B_1)[\mbox{\boldmath$b$}^1]$ of the outer point $B_1$ thus
\begin{equation}
\begin{gathered}
\mathbf{q}_2 \sim c \cdot \mathbf{b}_1 + \mathbf{p}_2 \in \mbox{\boldmath$b$}^1 \Leftrightarrow
c \cdot \mathbf{b}_1 \mbox{\boldmath$b$}^1+\mathbf{p}_2 \mbox{\boldmath$b$}^1=0 \Leftrightarrow
c=-\frac{\mathbf{p}_2 \mbox{\boldmath$b$}^1}{\mathbf{b}_1 \mbox{\boldmath$b$}^1} \Leftrightarrow \\
\mathbf{q}_2 \sim -\frac{\mathbf{p}_2 \mbox{\boldmath$b$}^1}{\mathbf{b}_1 \mbox{\boldmath$b$}^1}
\mathbf{b}_1+\mathbf{p}_2 \sim \mathbf{p}_2 (\mathbf{b}_1 \mbox{\boldmath$b$}^1) - \mathbf{b}_1 (\mathbf{p}_2 \mbox{\boldmath$b$}^1)=
\mathbf{p}_2 h_{33}-\mathbf{b}_1 h_{23},
\end{gathered} \tag{3.2}
\end{equation}
where $h_{ij}$ ($i,j=0,1,2,3)$) is the inverse of the corresponding Coxeter-Schl\"afli matrix (see (2.1), where $q=3, r=4$) of the orthoscheme $\mathcal{S}$. 
The hyperbolic distance $h(p)$ can be calculated by the following formula:
\[
\begin{gathered}
\cosh{h(p)}=\cosh{P_2Q_2}=\frac{- \langle {\mathbf{q}}_2, {\mathbf{p}}_2 \rangle }
{\sqrt{\langle {\mathbf{q}}_2, {\mathbf{q}}_2 \rangle \langle {\mathbf{p}}_2, {\mathbf{p}}_2 \rangle}}= \\ =\frac{h_{23}^2-h_{22}h_{33}}
{\sqrt{h_{22}\langle \mathbf{q}_2, \mathbf{q}_2 \rangle}} =
\sqrt{\frac{h_{22}~h_{33}-h_{23}^2}
{h_{22}~h_{33}}}.
\end{gathered} \tag{3.3}
\]
We get that the volume $Vol(\mathcal{O}^r(p))$, the maximal height $h(p)$ of congruent hyperballs and the
$\sum_{i=1}^6 Vol(\mathcal{H}^{h(p)}_i \cap \mathcal{O}^r))$ depend only on the parameter $p$ of the truncated regular tetrahedron $\mathcal{O}^r(p)$.
Moreover,
the volume of the hyperball pieces can be computed by the formula (2.1) and the volume of $\mathcal{O}^r(p)$ can be determined by the Theorem 2.2.
Therefore, the density $\delta(\mathcal{O}^r(h(p)))$ is depended only on parameter $p$. 

The domain of the density function $\delta(\mathcal{O}^r(h(p)))$ can be extended for $4<p \in \mathbb{R}$. The octahedron $\mathcal{O}^r(p)$ is realized in $\HYP$ for any given $p$
parameter ($4<p\in\mathbb{R}$) but the corresponding octahedron tiling exists only for parameters $5 \le p\in \mathbb{N}$. 

Finally, we obtain the graph of the smooth density function $\delta(\mathcal{O}^r(h(p)))$ and we get after its careful analysis (cf. Fig. 3) the following
\begin{theorem}
The density function $\delta(\mathcal{O}^r(h(p)))$, ($p\in (4,\infty)$)
attains its maximum at $p^{opt} \approx 4.11320$. It is strictly increasing on the interval $(4,p^{opt})$ and strictly decreasing on the interval $(p^{opt},\infty)$. 
Moreover, the optimal density
$\delta^{opt}(\mathcal{O}^r(h(p^{opt}))) \approx 0.83173$,
however these hyperball packing configurations are only locally optimal and cannot be extended to the entirety
of the hyperbolic spaces $\mathbb{H}^3$ (see Fig.~3).
\end{theorem}
\begin{corollary}
The density function $\delta(\mathcal{O}^r(h(p)))$, ($\mathbb{N} \ni p \ge 5$) attains its maximum at the parameter
$p=5$ and the corresponding congruent hyperball packing $\mathcal{B}^{h(5)}$ related to the regular truncated octahedra
can be extended to the entire hyperbolic space. The maximal density is $\delta(\mathcal{O}^r(h(5))) \approx 0.76893$.
\end{corollary}
\begin{figure}[ht]
\centering
\includegraphics[width=12cm]{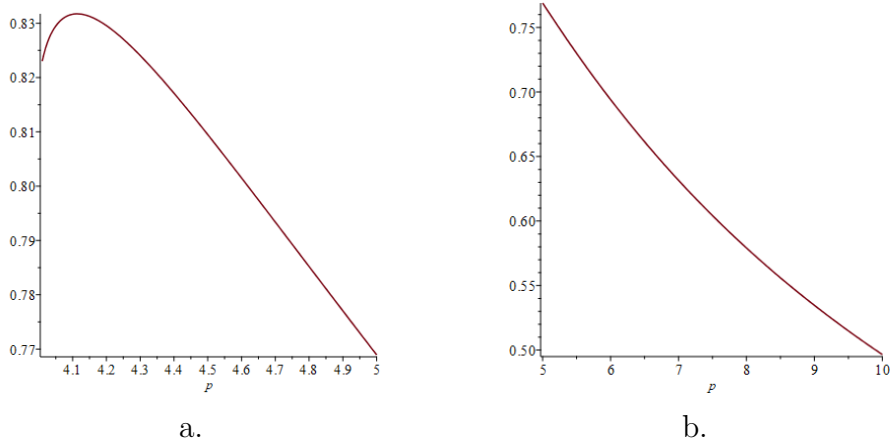}

a. \hspace{6cm} b.
\caption{a.~Density function $\delta(\mathcal{O}^r(h(p)))$, $(p\in [4,5])$.~b.~Density function $\delta(\mathcal{O}^r(h(p)))$, $(p\in[5,10]$).}
\label{}
\end{figure}
\begin{rmrk}We note here that these coincide with the hyperball packings to the regular prism tilings in $\HYP$ with Schl\"afli symbols
$\{p,4,3\}$ which are discussed in \cite{Sz06-1}.
\end{rmrk}
In the following Table we summarize the data of the hyperball packings for some parameters $p$, ($5 \ge p\in \mathbb{N}$).
\medbreak
\smallbreak
\centerline{\vbox{
\halign{\strut\vrule~\hfil $#$ \hfil~\vrule
&\quad \hfil $#$ \hfil~\vrule
&\quad \hfil $#$ \hfil\quad\vrule
&\quad \hfil $#$ \hfil\quad\vrule
&\quad \hfil $#$ \hfil\quad\vrule
\cr
\noalign{\hrule}
\multispan5{\strut\vrule\hfill\bf Table 1, $\{3,4,p\}$  \hfill\vrule}%
\cr
\noalign{\hrule}
\noalign{\vskip2pt}
\noalign{\hrule}
p & h(p) & Vol(\mathcal{O}^r(p))/48 & Vol(\mathcal{H}^{h(p)} \cap \mathcal{O}^r(p))/8 & \delta^{opt}(\mathcal{O}^r(h(p))) \cr
\noalign{\hrule}
5 & 0.69128565 & 0.16596371 & 0.12761435 & 0.76892924 \cr
\noalign{\hrule}
6 & 0.48121183 & 0.19616337 & 0.13616563 & 0.69414405 \cr
\noalign{\hrule}
7 & 0.37938071 & 0.21217704 & 0.13400462 & 0.63156984 \cr
\noalign{\hrule}
\vdots & \vdots  & \vdots  & \vdots  & \vdots \cr
\noalign{\hrule}
20 & 0.11318462 & 0.24655736 & 0.07142045 & 0.28967074 \cr
\noalign{\hrule}
\vdots & \vdots  & \vdots  & \vdots  & \vdots \cr
\noalign{\hrule}
50 & 0.04456095 & 0.25026133 & 0.03221956 & 0.12874366 \cr
\noalign{\hrule}
\vdots & \vdots  & \vdots  & \vdots  & \vdots \cr
\noalign{\hrule}
100 & 0.02223088 & 0.25078571 & 0.01676445 & 0.06684770 \cr
\noalign{\hrule}
p \to \infty & 0 & 0.25096025 & 0 & 0 \cr
\noalign{\hrule}}}}
\smallbreak
\subsection{Hyperball packings with non-congruent hyperballs in $3$-dimensional regular truncated octahedra}
We consider a regular truncated octahedron tiling $\mathcal{T}(\mathcal{O}^r(p))$ with Schl\"afli symbol $\{3,4,p\}$, $(5 \le p\in\mathbb{N})$.
One tile of it $\mathcal{O}^r(p)$ is illustrated in Fig.~2.a.b. 
This truncated octahedron can be derived also by truncation from a regular Euclidean octahedron centred at the origin with vertices $B_i$
($i\in\{1,\dots,6\})$. The truncating planes $\beta_i$ are the polar planes of outer vertices $B_i$ that can be the ultraparallel base planes 
of hyperballs $\mathcal{H}^{h_i(p)}_i$ $(i \in \{1,\dots,6\})$ with heights $h_i(p)$.

The distances between two base
planes $d(\beta_i,\beta_j)=:e_{ij}$ are equal ($i < j$, $i,j \in \{1,\dots 6\})$.
Moreover, the volume of the truncated simplex $\mathcal{O}^r(p)$ is denoted by $Vol(\mathcal{O}^r(p))$, similarly to the above section.

The distances of the plane $\beta_i$ ($i \in \{1,\dots 6\})$ from rectangular hexagon faces of the octahedron $\mathcal{O}^r(p)$ 
whose planes do not contain the vertex $B_i$ are equal and this distance is denoted by $w(p)$ (see Fig.~2.a.b). 
We construct non-congruent hyperball packings to $\mathcal{T}(\mathcal{O}^r(p))$ tilings therefore the hyperballs have to satisfy the
following requirements:
\begin{enumerate}
\item The base plane $\beta_i$ of the hyperball $\mathcal{H}^{h_i(p)}_i$ (with height $h_i(p)$) is the polar plane of the vertex $B_i$ (see Fig.~2),
\item $ card \{ int (\mathcal{H}^{h_i(p)}_i)\cap int(\mathcal{H}^{h_j(p)}_j)\}=0$, $i\ne j$,
\item $ card \{ int (\mathcal{H}^{h_i}_i \cap (B_jB_kB_l ~\text{plane})\}=0$ $(i,j,k,l \in \{1,\dots,6\},~ i~\ne j,k,l,~ j<k<l$ i.e.
the distance $w(p) \ge h_i(p)$.
\end{enumerate}
{\it If the hyperballs hold the above requirements then we obtain congruent or non-congruent hyperball packings $\mathcal{B}(\mathcal{O}^r(p))$ in hyperbolic $3$-space
derived by the structure of the considered Coxeter octahedron tilings.}

We introduce the locally density function $\delta(\mathcal{O}^r(p))$ related to $\mathcal{O}^r(p)$:
\begin{definition}
\begin{equation}
\delta(\mathcal{O}^r(p)):=\frac{\sum_{i=1}^6 Vol(\mathcal{H}^{h_i(p)}_i \cap \mathcal{O}^r(p))}{Vol({\mathcal{O}^r(p)})}. \notag
\end{equation}
\end{definition}
It is well known that a packing is locally optimal (i.e. its density is locally maximal), then it is locally stable i.e. each ball is fixed by
the other ones so that no ball of packing
can be moved alone without overlapping another ball of the given ball packing or by other requirements of the corresponding tiling.

To get the locally optimal non-congruent hyperball packing related to the regular octahedron tilings 
we distinguish two essential cases:
\begin{enumerate}
\item 
We set up from the optimal congruent ball arrangement $\mathcal{B}^{h(p)}$ (see former section)
where the neighbouring congruent hyperballs touch each other at the
''midpoints" of the edges of $\mathcal{O}^r(p)$.

We choose two opposite hyperballs (e.g. $\mathcal{H}^{h(p)}_5$ and $\mathcal{H}^{h(p)}_6$) and blow up
these hyperballs (hyperspheres) keeping the hyperballs $\mathcal{H}^{h_i(p)}_i$ $(i=1,2,3,4)$ tangent to them upto these hyperspheres touch each other and the symmetry 
plane $B_1B_2B_3B_4$ at $P_0$ ($t(p)=d(Q_0,P_0)$ see Fig.~2). During this expansion the heights of hyperballs $\mathcal{H}^{h_j(p)}_j(p)~(j=5,6)$ are 
$h_j(p)=h(p)+x$ where $x \in [0,\min\{h(p),w(p)-h(p), t(p)-h(p) \}]$. 
The height of further hyperballs are $h_1(p)=h_2(p)=h_3(p)=h_4(p)=h(p)-x$. (If $x=0$ then the hyperballs are congruent.) 

We extend this procedure to images of the hyperballs $\mathcal{H}^{h_i(p)}_i(p)$ $(i=1,\dots, 6)$ by the considered Coxeter group and obtain
non-congruent hyperball arrangements $\mathcal{B}^x_1(p)$.

Applying the Definition 3.5 we obtain the density function $\delta_1(\mathcal{O}^r(x,p))$:
\begin{equation}
\begin{gathered}
\delta_1(\mathcal{O}^r(x,p))=\frac{2 \cdot Vol(\mathcal{H}^{h(p)+x} \cap \mathcal{O}^r(p))+4 \cdot Vol(\mathcal{H}^{h(p)-x}
\cap \mathcal{O}^r(p))}{Vol({\mathcal{O}^r(p)})}, \\
\text{where} ~ x \in [0,\min\{h(p),w(p)-h(p), t(p)-h(p) \}].
\end{gathered} \tag{3.4}
\end{equation}
\item Now, we start from the non-congruent ball arrangement 
where two ''larger hyperballs" with base planes $\beta_5$ and $\beta_6$ are tangent at the centre $P_0$ of
the octahedron, while hyperballs at the remaining four vertices touch both ''larger"
hyperballs. The point of tangency of the above two larger hyperballs on line $B_5B_6$
is denoted by $P_0=I(0)$. We blow up
the hyperball $\mathcal{H}^{t(p)}_5(p)$ with base plane $\beta_5$ keeping the hyperballs $\mathcal{H}^{h_i(p)}_i(p)$ $(i=1,2,3,4,6)$ tangent to it. 
The point of tangency of the hyperballs $\mathcal{H}^{h_5(p)}_5(p)$ and $\mathcal{H}^{h_6(p)}_6(p)$
along line $B_5B_6$ is denoted by $I(x)$ where $x$ is the
hyperbolic distance between $P_0=I(0)$ and $I(x)$. During this expansion the height of hyperball $\mathcal{H}^{h_5(p)}(p)$ is 
$h_5(p)=t(p)+x$ where $x \in [0, 2h(p)-t(p)]$. 
We note here, that $w(p)\ge 2h(p)$ and $2t(p) \ge 2h(p)$, therefore the hyperball $\mathcal{H}^{h_5}_5(p)$ can be blown upto this hypersphere touches the planes 
$\beta_i$ $(i=1,2,3,4)$. 
\end{enumerate}
We extend this procedure to images of the hyperballs $\mathcal{H}^{h_i(p)}_i(p)$ $(i=1,\dots, 6)$ by the considered Coxeter group and obtain
non-congruent hyperball arrangements $\mathcal{B}^x_2(p)$.

Applying the Definition 3.5 we obtain the density function $\delta_2(\mathcal{O}^r(x,p))$:
{\footnotesize \begin{equation}
\begin{gathered}
\delta_2(\mathcal{O}^r(x,p))=\\ =\frac{Vol(\mathcal{H}^{t(p)+x} \cap \mathcal{O}^r(p))+Vol(\mathcal{H}^{t(p)-x} \cap \mathcal{O}^r(p))+4 \cdot Vol(\mathcal{H}^{2h(p)-t(p)-x}
\cap \mathcal{O}^r(p))}{Vol({\mathcal{O}^r(p)})}, \\
\text{where} ~ x \in [0, 2h(p)-t(p)].
\end{gathered} \tag{3.5}
\end{equation}
}
{\it The main problem is: what is the maximum of density functions $\delta_i(\mathcal{O}^r(x,p))$ $(i=1,2)$ for given integer parameters $p \ge 5$ where
$x \in \mathbb{R}$, and $x \in [0,\min\{h(p),w(p)-h(p), t(p)-h(p) \}]$ or $x \in [0, 2h(p)-t(p)]$}.

During this expansion process we can compute for a given integer parameter $p \ge 5$ the densities
$\delta_i(\mathcal{O}^r(x,p))$ of considered packings as the function of $x$.
\subsubsection{Computations for parameter $p \ge 5$}
Every $3$-dimensional hyperbolic truncated regular octahedron can be derived from a $3$-dimensional regular Euclidean octahedron.
We introduce a projective coordinate system (see Section 2 and Fig.~2.a) and a unit sphere $\mathbb{S}^{2}$ centred at the origin which is
interpreted as the ideal boundary of $\overline{\mathbb{H}}^3$ in Beltrami-Cayley-Klein's ball model.

Now, we consider a $3$-dimensional regular Euclidean octahedron centred at the origin with outer vertices regarding the Beltrami-Cayley-Klein's ball model.
The projective coordinates of the vertices of this tetrahedron are
\begin{equation}
\begin{gathered}
B_1=(1,y,0,0,0); ~ B_2=(1,0,y,0);~
B_3=(1,-y,0,0);~B_4=(1,0,-y,0); \\  B_5=(1,0,0,y);~B_6=(1,0,0,-y) ~\text{where} ~ ~ 1 < y\in \mathbb{R}.
\end{gathered} \tag{3.6}
\end{equation}
The truncated octahedron $\mathcal{O}^r(p)$ can be derived from the above regular octahedron by cuttings with the polar planes of 
vertices $B_i$ $(i=1,\dots,6)$. The images of $\mathcal{O}^r(p)$ under reflections on its side facets generate a tiling in $\HYP$
if its non-right dihedral angles are $\frac{2\pi}{p}$ $(5 \le p \in \mathbb{N})$. It is easy to see, that if the parameter $p$ is given, then
\begin{equation}
y=\frac{(3\cos{\frac{2\pi}{p}}+1)(\cos{\frac{2\pi}{p}}+1)}{{\cos{\frac{2\pi}{p}}+1}}. \tag{3.7}
\end{equation}
We have to determine for any parameter $p$ the distances $h(p)$ and $w(p)$. The values of $h(p)$ can be derived
from formula (3.3) and $w(p)$ follows from the next formula:
\begin{equation}
\begin{gathered}
\sinh{w(p)}=\Bigg|{\frac{\langle \bb_1,\bt\rangle}{\sqrt{-\langle \bb_1,\bb_1\rangle \langle \bt,\bt\rangle}}}\Bigg|=
\sqrt{\frac{y^4+3}{(3-y^2){(y^2-1)}}}. \tag{3.8}
\end{gathered}
\end{equation}
\begin{enumerate}
\item If $p=5$ then we obtain the following results:
\begin{equation}
\begin{gathered}
2h(5)\approx 1.38257; w(5) \approx 1.71082 \Rightarrow ~w(5) > 2h(5); ~ t(5) \approx 1.16974.
\end{gathered} \notag
\end{equation}
We note here, that if $x=0$ then the hyperspheres are congruent (see the former section).
Therefore, we can compute during the expansion process for the given parameter $p =5$ the densities of
$\delta_i(\mathcal{O}^r(x,5))$ $(i=1,2)$ (see (3.4-5)) of considered packings as the function of $x$ using the formulas (2.1), (3.4-8) and Theorem 2.2:
\begin{equation}
\begin{gathered}
\delta_1(\mathcal{O}^r(x,5))=\frac{2Vol(\mathcal{H}^{h(5)+x} \cap \mathcal{O}^r(5))+4 \cdot Vol(\mathcal{H}^{h(7)-x}
\cap \mathcal{O}^r(5))}{Vol({\mathcal{O}^r(5)})}, \\
\text{where} ~ x \in [0, t(5)-h(5) \approx 0.47845].
\end{gathered} \tag{3.9}
\end{equation}
{\begin{equation}
\begin{gathered}
\delta_2(\mathcal{O}^r(x,5))=\\ =(Vol(\mathcal{H}^{t(5)+x} \cap \mathcal{O}^r(5))+Vol(\mathcal{H}^{t(5)-x} \cap \mathcal{O}^r(5))+ \\ +4Vol(\mathcal{H}^{2h(5)-t(5)-x}
\cap \mathcal{O}^r(5)))/ {Vol({\mathcal{O}^r(5)})}, \\
\text{where} ~ x \in [0, 2h(5)-t(5)\approx 0.21285].
\end{gathered} \tag{3.10}
\end{equation}
}
The graphs of $\delta_i(\mathcal{O}^r(x,5))$ are described in Fig.~4. Analyzing the above density functions we get that their 
maximal densities are achieved at the starting point of the corresponding intervals. The maximal density belongs to the function $\delta_1(\mathcal{O}^r(x,5))$ 
with parameter $x=0$, $\delta_1(\mathcal{O}^r(x,5)) \approx 0.76893$ (see also Table 1). We note here, that the density of $\delta_1(\mathcal{O}^r(x,5))$ at the endpoint of the 
interval $[0, t(5)-h(5) \approx 0.47845]$ is
$\delta_1(\mathcal{O}^r(t(5)-h(5),5)) \approx 0.72624$ (see Fig.~4.a).
\begin{figure}[ht]
\centering
\includegraphics[width=12cm]{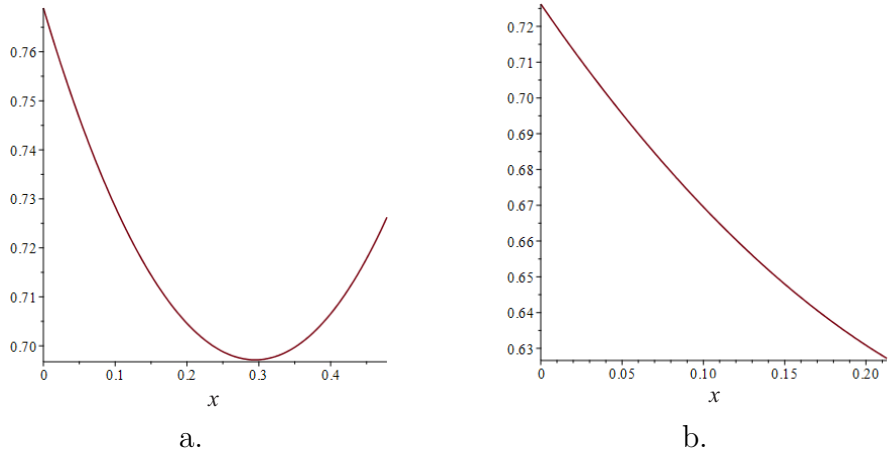}

a. \hspace{6cm} b.
\caption{a.~Density function $\delta_1(\mathcal{O}^r(x,5))$ where $x \in [0, t(5)-h(5) \approx 0.47845]$.
b.~Density function $\delta_2(\mathcal{O}^r(x,5))$ where $x \in [0, 2h(5)-t(5) \approx 0.21285]$.}
\label{}
\end{figure}
\item If $p > 5$ then $2h(5) < t(p) < w(p)$, therefore In this case we have to only examine the density function $\delta_1(\mathcal{O}^r(x,p))$ (see (3.4)). 
Similarly to the above computations for parameter $p=5$ we can analyze the density function and their maximum of
non-congruent hyperball packings generated by considered truncated octahedron tilings (or Coxeter tilings $\{3,4,p\}$) for all possible integer parameters
$p > 5$. In this case we have to examine the density function (3.4). Using the results of Theorem 3.2 and Corollary 3.3 we obtain the following
\end{enumerate}
\begin{theorem}
\begin{enumerate}
\item The maximum of the density function $\delta_1(\mathcal{O}^r(x,p))$ is attained at the starting
point of the corresponding intervals $x \in [0,h(p)]$ or $x \in [0,t(p)-h(p)]$ depending on integer parameter $p\ge 5$ i.e.
the congruent hyperball packing provides the densest hyperball packing for a given parameter $p$.
\item The maximum of the density function $\delta_2(\mathcal{O}^r(x,5))$ is achieved at the starting point of the interval $[0, 2h(5)-t(5) \approx 0.21285]$.
If $p>5$ then this case does not occur because $2h(5) < t(p) < w(p)$.
\item The maximum of the density functions $\delta_i(\mathcal{O}^r(x,p))$ ($p \ge 5$, integer parameter) is achieved at the parameters $x=0$, $p=5$.
Therefore, the density upper bound of the congruent and non-congruent hyperball packings is $\approx 0.76893$.
\end{enumerate}
\end{theorem}
\begin{figure}[ht]
\centering
\includegraphics[width=14cm]{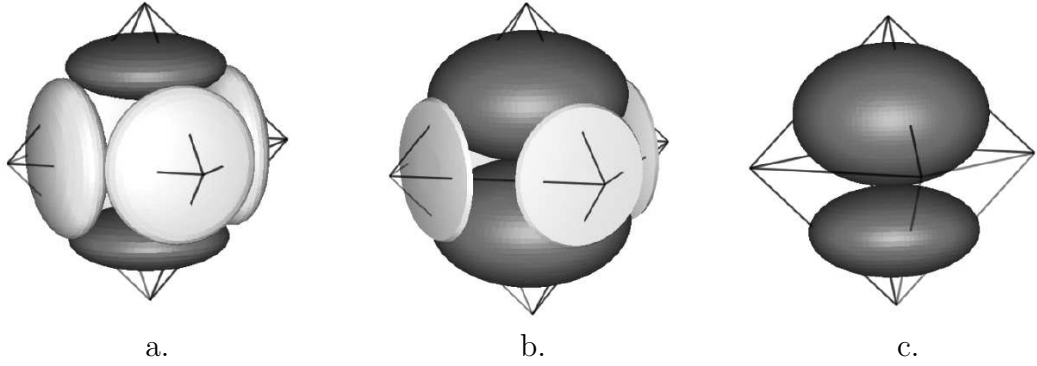}

\hspace{0.5cm} a. \hspace{4.3cm} b. \hspace{4.3cm} c.
\caption{a.~The densest packing configuration with density $\delta_1(\mathcal{O}^r(0,5))\approx 0.76893$. Here the hyperballs are congruent.
b.~The packing arrangement of parameters $p=5$, $x=t(5)-h(5) \approx 0.47845$ with density $\approx 0.72624$ where the two ''larger hyperballs" 
with base planes $\beta_5$ and $\beta_6$ are tangent at the centre $P_0$ of the octahedron.~c.~
The largest hypersphere $\mathcal{H}^{h_5(5)}_5(5)$ touches the planes $\beta_i$ $(i=1,2,3,4)$ and the opposite hypersphere $\mathcal{H}^{h_6(5)}_6(5)$.}
\label{}
\end{figure}
\subsection{Hyperball packings with congruent hyperballs related to regular truncated cube tilings $\{4,3,p\}$}
We consider a regular truncated cube tiling $\mathcal{T}(\mathcal{C}^r(p))$ with Schl\"afli symbol $\{4,3,p\}$, $(7 \le p \in \mathbb{N})$.
These cube tilings are derived by duality from the Coxeter tilings $\{p,3,4\}$ whose fundamental 
domains are simply truncated orthoschems (e.g. $P_1P_2P_3Q_1Q_2Q_3$ in Fig.~6.b).
One tile of it $\mathcal{C}^r(p)$ (a truncated cube) is illustrated in Fig.~6.a.b which can be derived also by truncation from a regular
Euclidean cube centred at the origin with vertices $B_i$
($i\in\{1,\dots,8\})$. The truncating planes $\beta_i$ are the polar planes of outer vertices $B_i$ that can be the ultraparallel base planes 
of hyperballs $\mathcal{H}^{s}_i$ $(i \in \{1,\dots,8\})$ with height $s$. The non-orthogonal dihedral angles of $\mathcal{C}^r(p)$ are equal to $\frac{2\pi}{p}$.
The distances between two base planes are equal $d(\beta_i,\beta_j)=:e_{ij}=2h(p)$ ($i < j$, $i,j \in \{1,\dots 8\})$ ($d$ is the hyperbolic distance function)
therefore the height of a hyperball is at most $h(p)$ (see Fig.~6.a.b). 
It is clear, that in the congruent, densest case the heights of the hyperballs are $h(p)$ i.e. the neighbouring hyperballs touch each other.

We consider a saturated congruent hyperball packing $\mathcal{B}^{h(p)}$ of hyperballs $\mathcal{H}^{h(p)}_i$ related to the above regular, truncated cube tilings 
$\{4,3,p\}$ $(7 \le p \in \mathbb{N})$.
\begin{figure}[ht]
\centering
\includegraphics[width=6.5cm]{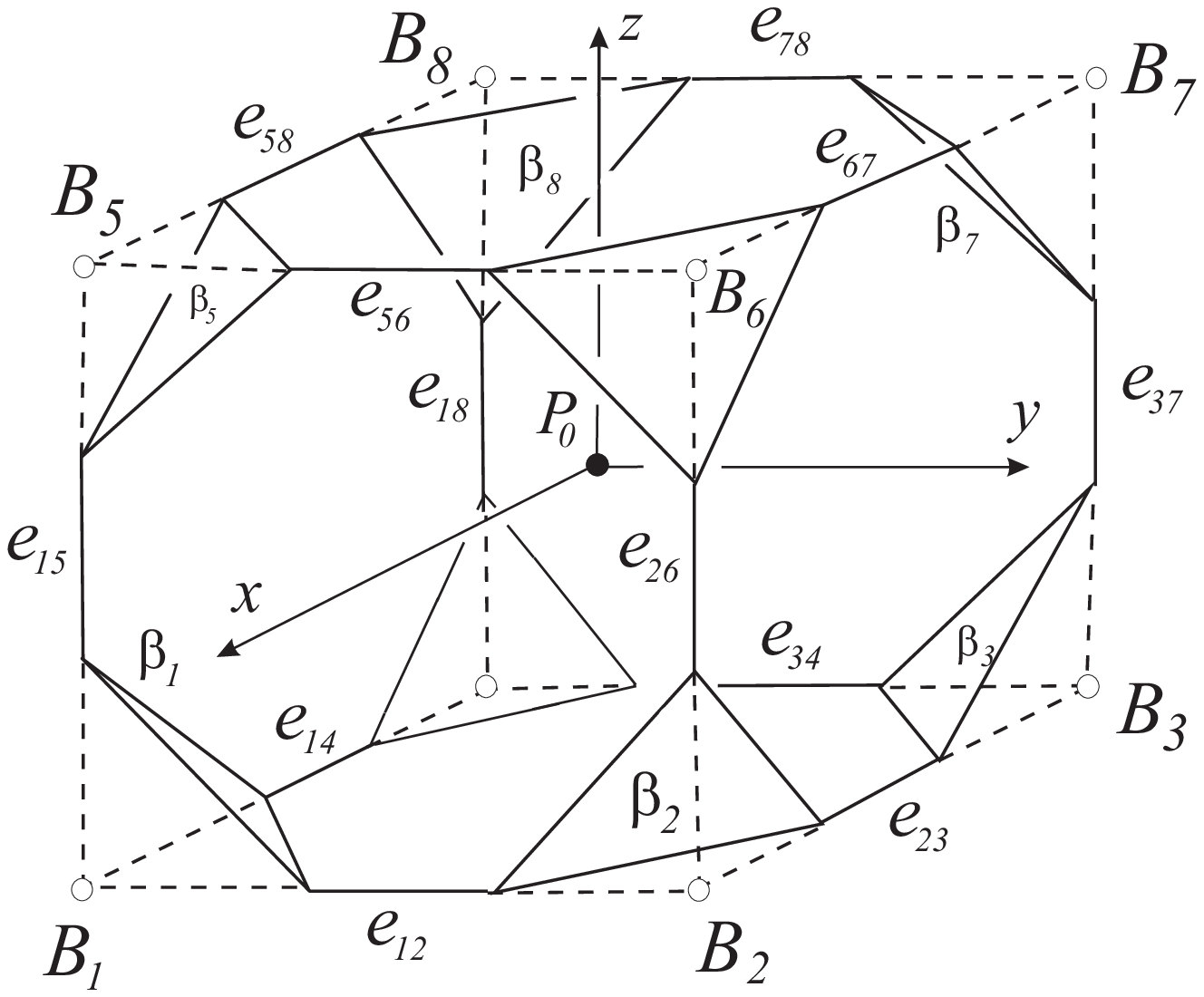} \includegraphics[width=6.5cm]{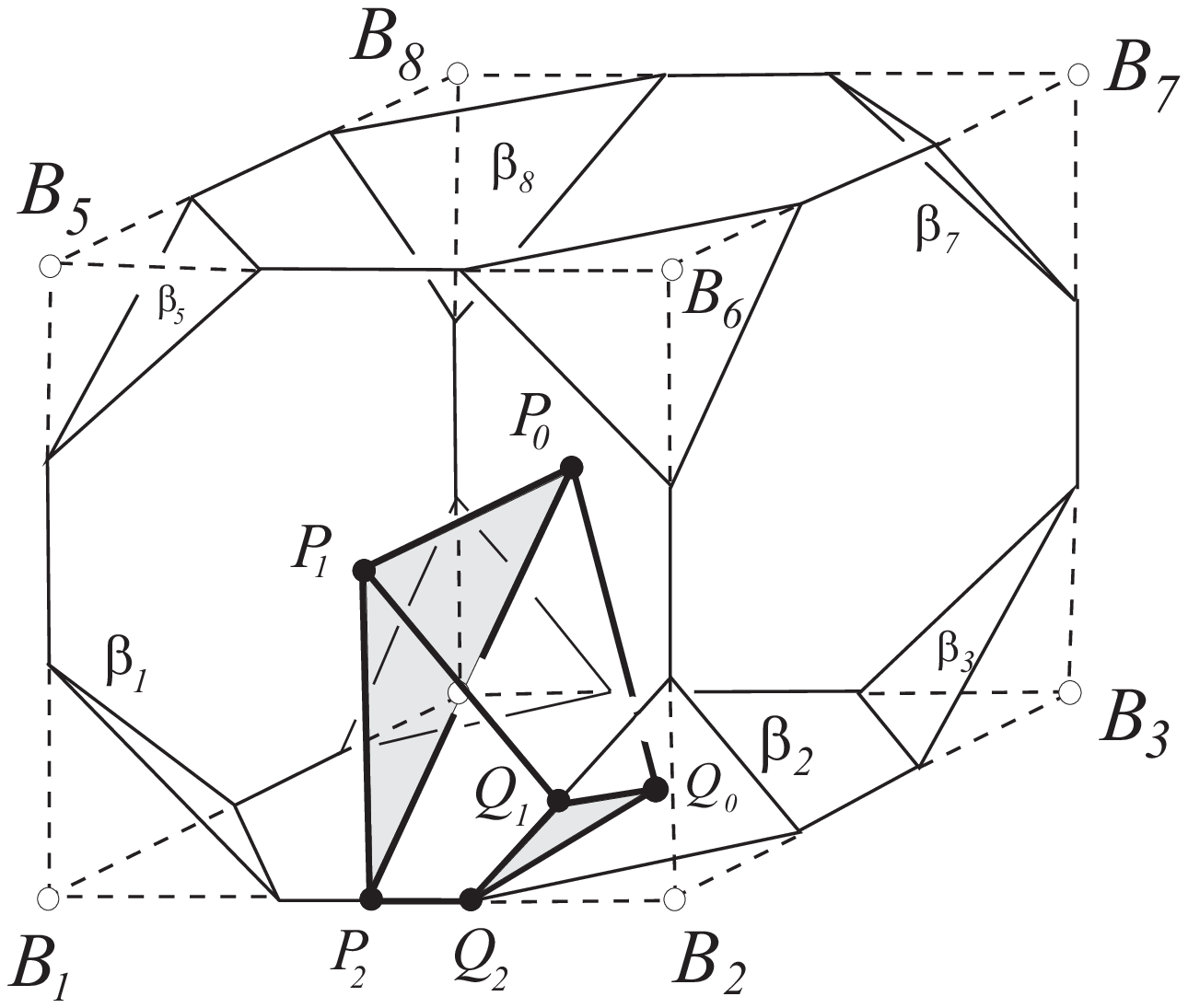}

a. \hspace{6cm} b.
\caption{Regular truncated cube. ~b.~Regular truncated cube with complete orthoscheme of degree $m=1$ (simple frustum orthoscheme)
with outer vertex $B_2$. This orthoscheme is truncated by its polar plane $\beta_2=pol(B_2)$.}
\label{}
\end{figure}
The volume of the truncated cube $\mathcal{C}^r(p)$ is denoted by $Vol(\mathcal{C}^r(p))$ and
we introduce the locally density function $\delta(\mathcal{C}^r(h(p)))$ related to $\mathcal{C}^r(p)$:
\begin{definition}
\begin{equation}
\delta(\mathcal{C}^r(h(p))):=\frac{\sum_{i=1}^8 Vol(\mathcal{H}_i^{h(p)} \cap \mathcal{C}^r(p))}{Vol({\mathcal{C}^r(p)})}=
\frac{8 Vol(\mathcal{H}_i^{h(p)} \cap \mathcal{C}^r(p))}{Vol({\mathcal{C}(p)})}. \notag
\end{equation}
\end{definition}
If the parameter $p$ is given then the common length of the common perpendiculars $2h(p)=e_{ij}$ $(i < j$, $i,j \in \{1,\dots,8\})$
can be determined by the machinery of the projective geometry, similarly to the octahedral case.
\[
\begin{gathered}
\cosh{h(p)}=\cosh{P_2Q_2}=\frac{- \langle {\mathbf{q}}_2, {\mathbf{p}}_2 \rangle }
{\sqrt{\langle {\mathbf{q}}_2, {\mathbf{q}}_2 \rangle \langle {\mathbf{p}}_2, {\mathbf{p}}_2 \rangle}}= \\ =\frac{h_{23}^2-h_{22}h_{33}}
{\sqrt{h_{22}\langle \mathbf{q}_2, \mathbf{q}_2 \rangle}} =
\sqrt{\frac{h_{22}~h_{33}-h_{23}^2}
{h_{22}~h_{33}}}.
\end{gathered} \tag{3.11}
\]
where $h_{ij}$ ($i,j=0,1,2,3)$) is the inverse of the corresponding Coxeter-Schl\"afli matrix (see (2.2), where $q=3, r=4$) of the orthoscheme $P_1P_2P_3Q_1Q_2Q_3$ (see Fig.~6.b). 

The volume $Vol(\mathcal{C}^r(p))$ can be calculated by Theorem 2.2 and the volume of the hyperball pieces lying in $\mathcal{C}^r(p)$ can be computed by the formula (2.1) 
for each given parameter $p$ therefore the maximal height $h(p)$ of the congruent hyperballs and the
$\sum_{i=1}^8 Vol(\mathcal{H}_i^{h(p)} \cap \mathcal{C}^r(p)))$ depend only on the parameter $p$ of the truncated regular cube $\mathcal{C}^r(p)$.

Thus, the density $\delta(\mathcal{C}^r(h(p)))$ is depended only on parameter $p$ $(7 \le p\in \mathbb{N})$. 

Finally, we obtain after careful analysis of the smooth
density function (cf. Fig. 7) the following
\begin{theorem}
The density function $\delta(\mathcal{C}^r(h(p)))$, ($p\in (6,\infty)$)
attains its maximum at $p^{opt} \approx 6.33962$. It is strictly increasing on the interval $(6,p^{opt})$ and strictly decreasing on the interval 
$(p^{opt},\infty)$. Moreover, the optimal density
$\delta^{opt}(\mathcal{C}^r(h(p^{opt}))) \approx 0.70427$,
however these hyperball packing configurations are only locally optimal and cannot be extended to the entirety
of the hyperbolic spaces $\mathbb{H}^3$ (see Fig.~7).
\end{theorem}
\begin{corollary}
The density function $\delta(\mathcal{C}^r(h(p)))$, ($7 \le p \in \mathbb{N}$) attains its maximum at the parameter
$p=7$. The congruent hyperball packing $\mathcal{B}^{h(7)}$ related to the regular truncated cube tilings
can be extended to the entire hyperbolic space. The maximal density is $\delta(\mathcal{C}^r(h(7))) \approx 0.68839$.
\end{corollary}
\begin{figure}[ht]
\centering
\includegraphics[width=12cm]{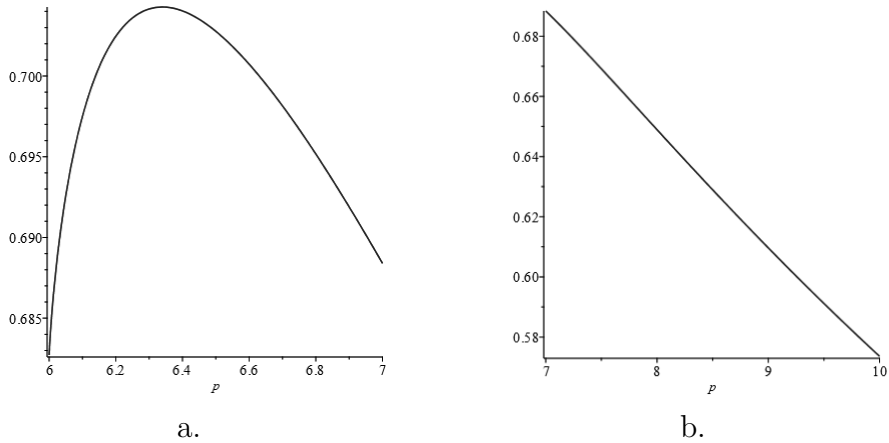}

a. \hspace{6cm} b.
\caption{a.~Density function $\delta(\mathcal{C}^r(h(p)))$, $(p\in [6,7])$.~b.~Density function $\delta(\mathcal{C}^r(h(p)))$, $(p\in[7,10]$).}
\label{}
\end{figure}
\begin{rmrk}We note here that these results coincide with the hyperball packings to the regular prism tilings in $\HYP$ with Schl\"afli symbols
$\{p,3,4\}$ which are discussed in \cite{Sz06-1}.
\end{rmrk}
In the following Table we summarize the data of the hyperball packings for some parameters $p$, ($5 \ge p\in \mathbb{N}$).
\medbreak
\smallbreak
\centerline{\vbox{
\halign{\strut\vrule~\hfil $#$ \hfil~\vrule
&\quad \hfil $#$ \hfil~\vrule
&\quad \hfil $#$ \hfil\quad\vrule
&\quad \hfil $#$ \hfil\quad\vrule
&\quad \hfil $#$ \hfil\quad\vrule
\cr
\noalign{\hrule}
\multispan5{\strut\vrule\hfill\bf Table 2, $\{4,3,p\}$  \hfill\vrule}%
\cr
\noalign{\hrule}
\noalign{\vskip2pt}
\noalign{\hrule}
p & h(p) & Vol(\mathcal{C}^r)/48 & Vol(\mathcal{H}^{h(p)} \cap \mathcal{C}^r(p))/6 & \delta^{opt}(\mathcal{C}^r(h(p))) \cr
\noalign{\hrule}
7 & 1.03799291 & 0.16297337 & 0.11218983 & 0.68839367 \cr
\noalign{\hrule}
8 & 0.76428546 & 0.18789693 & 0.12193107 & 0.64892530 \cr
\noalign{\hrule}
9 & 0.62216938 & 0.20295023 & 0.12372607 & 0.60963750 \cr
\noalign{\hrule}
\vdots & \vdots  & \vdots  & \vdots  & \vdots \cr
\noalign{\hrule}
20 & 0.23086908 & 0.24206876 & 0.08613744 & 0.35583872 \cr
\noalign{\hrule}
\vdots & \vdots  & \vdots  & \vdots  & \vdots \cr
\noalign{\hrule}
50 & 0.08938872 & 0.24956032 & 0.04129724 & 0.16547999 \cr
\noalign{\hrule}
\vdots & \vdots  & \vdots  & \vdots  & \vdots \cr
\noalign{\hrule}
100 & 0.04449475 & 0.25061105 & 0.02191401 & 0.08744233 \cr
\noalign{\hrule}
p \to \infty & 0 & 0.25096025 & 0 & 0 \cr
\noalign{\hrule}}}}
\subsection{Hyperball packings with non-congruent hyperballs in $3$-dimensional regular truncated cubes}
We consider a regular truncated cube tiling $\mathcal{T}(\mathcal{C}^r(p))$ with Schl\"afli symbol $\{4,3,p\}$, $(7 \le p\in\mathbb{N})$.
One tile of it $\mathcal{C}^r(p)$ is illustrated in Fig.~6.a.b. 
This truncated cube can be derived also by truncation from a regular Euclidean cube centred at the origin with vertices $B_i$
($i\in\{1,\dots,8\})$. The truncating planes $\beta_i$ are the polar planes of outer vertices $B_i$ that can be the ultraparallel base planes 
of hyperballs $\mathcal{H}^{h_i(p)}_i$ $(i \in \{1,\dots,8\})$ with heights $h_i(p)$.
The distances between two base
planes $d(\beta_i,\beta_j)=:e_{ij}=2h(p)$ are equal ($i < j$, $i,j \in \{1,\dots 8\})$ and can be determined by formula (3.11).
Moreover, the volume of the truncated cube $\mathcal{C}^r(p)$ is denoted by $Vol(\mathcal{C}^r(p))$, and can be computed by Theorem 2.2, similarly to the above section.

The distances of the plane $\beta_i$ ($i \in \{1,\dots 8\})$ from rectangular octagon faces of the cube $\mathcal{C}^r(p)$ 
whose planes do not contain the vertex $B_i$ are equal and this distance is denoted by $w(p)$ (see Fig.~6.a.b). 

We would like to construct non-congruent hyperball packings to $\mathcal{T}(\mathcal{C}^r(p))$ tilings therefore the hyperballs have to satisfy the
following requirements:
\begin{enumerate}
\item The base plane $\beta_i$ of the hyperball $\mathcal{H}^{h_i(p)}_i$ (with height $h_i(p)$) is the polar plane of the vertex $B_i$ (see Fig.~6),
\item $ card \{ int (\mathcal{H}^{h_i(p)}_i)\cap int(\mathcal{H}^{h_j(p)}_j)\}=0$, $i\ne j$,
\item $ card \{ int (\mathcal{H}^{h_i(p)}_i \cap B_jB_kB_l ~\text{plane}\}=0$ $(i,j,k,l \in \{1,\dots,8\},~ i~\ne j,k,l,~ j<k<l$ i.e.
the distance $w(p) \ge h_i(p)$.
\end{enumerate}
{\it If the hyperballs hold the above requirements then we obtain congruent or non-congruent hyperball packings in the cube $\mathcal{C}^r(p)$ 
and if we extend them by the structure of the considered Coxeter cube tilings $\mathcal{T}(\mathcal{C}^r(p))$ then we get 
hyperball packings $\mathcal{B}(p)$ in hyperbolic $3$-space.}

We introduce the locally density function $\delta(\mathcal{C}^r(p))$ related to above packings:
\begin{definition}
\begin{equation}
\delta(\mathcal{C}^r(p)):=\frac{\sum_{i=1}^8 Vol(\mathcal{H}^{h_i(p)}_i \cap \mathcal{C}^r(p))}{Vol({\mathcal{C}^r(p)})}. \notag
\end{equation}
\end{definition}
We will use that a packing is locally optimal (i.e. its density is locally maximal), then it is locally stable i.e. each ball is fixed by
the other ones so that no ball of packing
can be moved alone without overlapping another ball of the given ball packing or by other requirements of the corresponding tiling.

To get the locally optimal non-congruent hyperball packing arrangement related to the regular cube or cube tilings 
we distinguish three essential cases:
\begin{enumerate}
\item 
We set up from the optimal congruent ball arrangement (see former subsection)
where the neighbouring congruent hyperballs touch each other at the
''midpoints" of the edges of $\mathcal{C}^r(p)$.

We choose two opposite hyperballs (e.g. $\mathcal{H}^{h(p)}_2$ and $\mathcal{H}^{h(p)}_8$) and blow up
these hyperballs (hyperspheres) keeping the hyperballs $\mathcal{H}^{h_i(p)}_i(p)$ $(i\in\{1,3,4,5,6,7\})$ tangent to them upto their heigths $h_2(p)=h_8(p)=\min\{2h(p),w(p),t(p),s(p)\}$ 
where $s(p)=d(Q_1,P_1)$ and $t(p)=d(Q_0,P_0)$ (see Fig.~6).
During this expansion the heights of hyperballs $\mathcal{H}^{h_j(p)}_j~(j=2,8)$ are $h_j(p)=h(p)+x$ where $x\in [0, \min\{2h(p),w(p), t(p)\}-h(p)]$. 
The heights of further hyperballs are $h_1(p)=h_3(p)=h_4(p)=h_5(p)=h_6(p)=h_7(p)=h(p)-x$. (If $x=0$ then the hyperballs are congruent.) 

We extend this procedure to images of the hyperballs $\mathcal{H}^{h_i(p)}_i$ $(i \in \{1,\dots, 8\})$ by the considered Coxeter group and obtain
non-congruent hyperball arrangements $\mathcal{B}^x_1(p)$.

Applying the Definition 3.11 we obtain the density function $\delta_1(\mathcal{C}^r(x,p))$:
\begin{equation}
\begin{gathered}
\delta_1(\mathcal{C}^r(x,p))=\frac{2 \cdot Vol(\mathcal{H}^{h(p)+x} \cap \mathcal{C}^r(p))+6 \cdot Vol(\mathcal{H}^{h(p)-x}
\cap \mathcal{C}^r(p))}{Vol({\mathcal{C}^r(p)})}, \\
\text{where} ~ x \in [0,\min\{2h(p),w(p), t(p)\}-h(p)].
\end{gathered} \tag{3.12}
\end{equation}
\item Now, we start from the non-congruent ball arrangement 
where two ''larger hyperballs" with base planes $\beta_2$ and $\beta_8$ are tangent at the centre $P_0$ of
the cube, while hyperballs at the remaining six vertices touch the corresponding ''larger"
hyperball. The point of tangency of the above two larger hyperballs on line $B_2B_8$
is denoted by $P_0=I(0)$. We blow up
the hyperball (hypersphere) $\mathcal{H}^{t(p)}_2$ with base plane $\beta_2$ keeping the hyperballs $\mathcal{H}^{h_i(p)}_i$ $(i=1,3,6,8)$ tangent to it while 
the hyperballs $\mathcal{H}^{h_i(p)}_i$ $(i=4,5,7)$ are blowing up touching the hyperball $\mathcal{H}^{h_8(p)}_8$. 

The point of tangency of the hyperballs $\mathcal{H}^{h_2(p)}_2$ and $\mathcal{H}^{h_8(p)}_8$
along line $B_2B_8$ is denoted by $I(x)$ where $x$ is the
hyperbolic distance between $P_0=I(0)$ and $I(x)$. During this expansion the height of hyperball $\mathcal{H}^{h_2(p)}$ is 
$h_2(p)=t(p)+x$ where $x \in [0, \min\{ 2h(p)-t(0), t(p), w(p)-t(p),s(p)-h(p) \} $.  
We note here, that $w(p)\ge 2h(p)$ and $t(p) \ge s(p) \ge h(p)$, therefore the hyperball $\mathcal{H}^{h_2(p)}_2$ can be at most blown upto this hypersphere touch the planes 
$\beta_i$ $(i=1,2,3,4)$. 
We extend this procedure to images of the hyperballs $\mathcal{H}^{h_i(p)}_i$ $(i \in \{1,\dots, 8\})$ by the considered Coxeter group and obtain
non-congruent hyperball arrangements $\mathcal{B}^x_2(p)$.

Applying the Definition 3.11 we obtain the density function $\delta_2(\mathcal{C}^r(x,p))$:
{\begin{equation}
\begin{gathered}
\delta_2(\mathcal{C}^r(x,p))=(Vol(\mathcal{H}^{t(p)+x} \cap \mathcal{C}^r(p))+\\+Vol(\mathcal{H}^{t(p)-x} \cap \mathcal{C}^r(p))+3 \cdot Vol(\mathcal{H}^{2h(p)-t(p)-x}
\cap \mathcal{C}^r(p))+\\+3 \cdot Vol(\mathcal{H}^{2h(p)-t(p)+x}
\cap \mathcal{C}^r(p)))/Vol({\mathcal{C}^r(p)}), \\
\text{where} ~ x \in [0, \min\{ 2h(p)-t(0), t(p), w(p)-t(p),s(p)-h(p) \} .
\end{gathered} \tag{3.13}
\end{equation}
}
\item We set up from the congruent ball arrangement (see former section)
where the neighbouring congruent hyperballs touch each other at the
''midpoints" of the edges of $\mathcal{C}^r(p)$. We distinguish two different classes of hyperballs related to two complementary tetrahedral sublattices of vertices of the
cube. E.g. the hyperballs $\mathcal{H}^{h_i(p)}_i$ ($i\in \{1,3,6,8\}$) form the first class and the remaining hyperballs form the second class.
We blow up the hyperballs of the first class keeping the hyperballs $\mathcal{H}^{h_i(p)}_i$ $(i\in\{2,4,5,7\})$ tangent to them upto their heigths 
$h_1(p)=h_3(p)=h_4(p)=h_7(p) = \min\{2h(p),s(p)\}$ (see Fig.~6).
During this expansion the heights of hyperballs $\mathcal{H}^{h_j(p)}_j~(j\in \{1,3,6,8\}$ are $h_j(p)=h(p)+x$ and 
$\mathcal{H}^{h_j(p)}_j~(j\in \{2,4,5,7\})$ are $h_j(p)=h(p)-x$ where $x\in [0, \min\{h(p),s(p)-h(p)]$. 

We extend this procedure to images of the hyperballs $\mathcal{H}^{h_i(p)}_i$ $(i \in \{1,\dots, 8\})$ by the considered Coxeter group and obtain
non-congruent hyperball arrangements $\mathcal{B}^x_3(p)$.
Applying the Definition 3.5 we obtain the density function $\delta_2(\mathcal{C}^r(x,p))$:
\begin{equation}
\begin{gathered}
\delta_3(\mathcal{C}^r(x,p))=\frac{4 Vol(\mathcal{H}^{h(p)+x} \cap \mathcal{C}^r(p))+4 \cdot Vol(\mathcal{H}^{h(p)-x}
\cap \mathcal{C}^r(p))}{Vol({\mathcal{C}^r(p)})}, \\
\text{where} ~ x \in [0,\min\{h(p),s(p)-h(p)].
\end{gathered} \tag{3.14}
\end{equation}
\end{enumerate}
{\it The main problem is: what is the maximum of density functions $\delta_i(\mathcal{C}^r(x,p))$ $(i\in \{1,2,3\})$ for given integer parameters $p \ge 7$ where
$x\in \mathbb{R}$ belongs to the corresponding intervals.}

During this expansion processes we can compute for a given integer parameter $p \ge 7$ the densities
$\delta_i(\mathcal{C}^r(x,p))$ of considered packings as the function of $x$.
\subsubsection{Computations}
Every $3$-dimensional hyperbolic truncated regular cube can be derived from a $3$-dimensional regular Euclidean cube (see Fig.~6).
We introduce a projective coordinate system (see Section 2 and Fig.~6.a) and a unit sphere $\mathbb{S}^{2}$ centred at the origin which is
interpreted as the ideal boundary of $\overline{\mathbb{H}}^3$ in Beltrami-Cayley-Klein's ball model.

Now, we consider a $3$-dimensional regular Euclidean cube centred at the origin with outer vertices regarding the Beltrami-Cayley-Klein's ball model.
The projective coordinates of the vertices of this cube are
\begin{equation}
\begin{gathered}
B_1=(1,y,-y,-y); ~ B_2=(1,y,y,-y);~B_3=(1,-y,y,-y);\\ B_4=(1,-y,-y,-y);~ B_5=(1,y,-y,y); \\ 
B_6=(1,y,y,y);~B_7=(1,-y,y,y);~B_8=(1,-y,-y,y)~\text{where} ~ ~ \frac{1}{\sqrt{3}} < y \in \mathbb{R}.
\end{gathered} \tag{3.15}
\end{equation}
The truncated cube $\mathcal{C}^r(p)$ can be derived from the above regular cube by cuttings with the polar planes of 
vertices $B_i$ $(i\in \{1,\dots,8\})$. The images of $\mathcal{C}^r(p)$ under reflections on its side facets generate a tiling in $\HYP$
if its non-right dihedral angles are $\frac{2\pi}{p}$ $(\mathbb{N} \ni p\ge 7)$. It is easy to see, that if the parameter $p$ is given, then
\begin{equation}
y=\sqrt{\frac{\cos{\frac{2\pi}{p}}}{\cos{\frac{2\pi}{p}}+1}}. \tag{3.16}
\end{equation}
We have to determine for any parameter $p$ the distances $h(p)$, $t(p)$, $s(p)$ and $w(p)$. The values of $h(p)$ can be derived
from formula (3.11). The distances $t(p)$ and $s(p)$ can be determined similarly to (3.11). 
$w(p)$ follows from the next formula:
\begin{equation}
\begin{gathered}
\sinh{w(p)}=\Bigg|{\frac{\langle \bb_1,\bt\rangle}{\sqrt{-\langle \bb_1,\bb_1\rangle \langle \bt,\bt\rangle}}}\Bigg|=
\sqrt{\frac{3y^4+1}{(1-3y^2){(y^2-1)}}}. \tag{3.17}
\end{gathered}
\end{equation}
\begin{enumerate}
\item If $p=7$ then we obtain the following results:
\begin{equation}
\begin{gathered}
2h(7)\approx 2.07599; w(7) \approx 2.07599 \Rightarrow ~w(7) = 2h(7); \\ t(7) \approx 1.67069,~s(7) \approx 1.03799.
\end{gathered} \notag
\end{equation}
We can compute during the expansion processes the densities of
$\delta_i(\mathcal{C}^r(x,7))$ $(i\in\{1,2,3\})$ (see Definition 3.11) of considered packings as the function of $x$ using the formulas (2.1), (3.11-14), (3.16-17) 
and Theorem 2.2.
\begin{equation}
\begin{gathered}
\delta_1(\mathcal{C}^r(x,7))=\frac{2 \cdot Vol(\mathcal{H}^{h(7)+x} \cap \mathcal{C}^r(7))+6 \cdot Vol(\mathcal{H}^{h(7)-x}
\cap \mathcal{C}^r(7))}{Vol({\mathcal{C}^r(7)})}, \\
\text{where} ~ x \in [0,\min\{2h(7),w(7),t(7)\}-h(7)\}\approx 0.63270].
\end{gathered} \tag{3.18}
\end{equation}
\begin{figure}[ht]
\centering
\includegraphics[width=12cm]{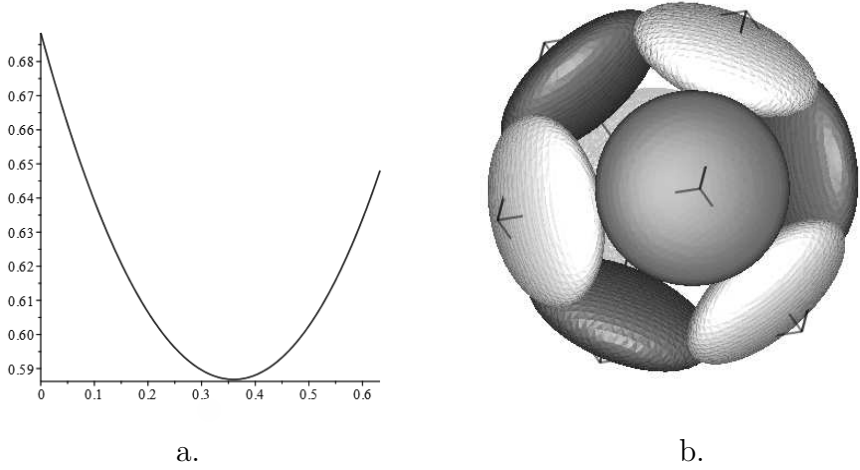}

a. \hspace{6cm} b.
\caption{a.~Density function $\delta_1(x,7))$ $(x\in [0, t(7)-h(7)],7)$. If $x=0$ then the density is $\approx 0.68839$. 
The density in the endpoint of the above interval is $\delta_2(\mathcal{C}^r(t(7)-h(7),7) \approx 0.64805$.~b.~The congruent hyperball arrangement with parameters
$x=0$, $p=7$ with density $\approx 0.68839$.}
\label{}
\end{figure}
The graph of $\delta_1(\mathcal{C}^r(x,7))$ are described in Fig.~8.a. Analyzing the above density function we get that the
maximal density is achieved at the starting point of the above interval (at the congruent case)
with density $\approx 0.68839$ (see Table 2). The density in the endpoint of the above interval is
$\delta_1(\mathcal{C}^r(t(7)-h(7)\approx 0.63270,7)) \approx 0.64805$.
{\begin{equation}
\begin{gathered}
\delta_2(\mathcal{C}^r(x,7))=(Vol(\mathcal{H}^{t(7)+x} \cap \mathcal{C}^r(7))+\\+Vol(\mathcal{H}^{t(7)-x} \cap \mathcal{C}^r(7))+3 \cdot Vol(\mathcal{H}^{2h(7)-t(7)-x}
\cap \mathcal{C}^r(7))+\\+3 \cdot Vol(\mathcal{H}^{2h(7)-t(7)+x}
\cap \mathcal{C}^r(7)))/Vol({\mathcal{C}^r(7)}), ~\text{where} \\ x \in [0, \min\{ 2h(7)-t(0), t(7), w(7)-t(7),s(7)-h(7) \}\approx 0.40530] .
\end{gathered} \tag{3.19}
\end{equation}
}
\begin{figure}[ht]
\centering
\includegraphics[width=12cm]{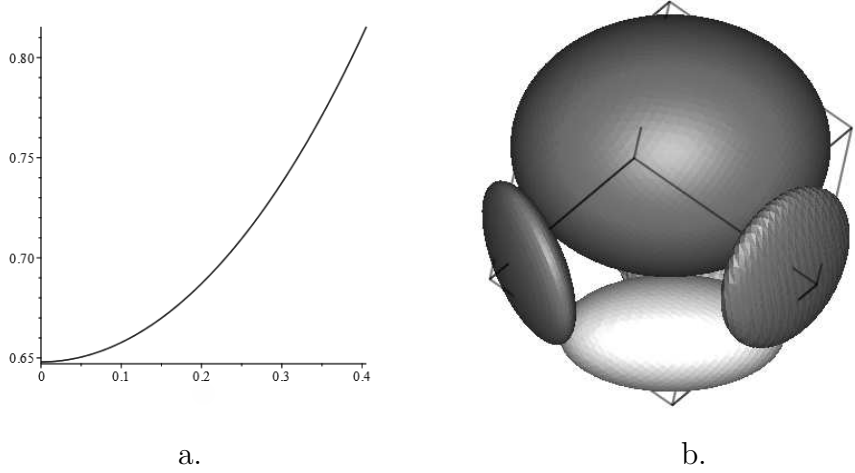}

a. \hspace{6cm} b.
\caption{a.~Density function $\delta_2(x,7))$ $(x\in [0, w(7)-t(7)],7)$. If $x=w(7)-t(7) \approx 0.40530$ then the density is $\approx 0.81542$. 
The density in the starting point of the above interval is $\delta_2(\mathcal{C}^r(0,7) \approx 0.64805$.~b.~The non-congruent hyperball arrangement with parameters
$x=w(7)-t(7) \approx 0.40530$, $p=7$ with density $\approx 0.81542$.}
\label{}
\end{figure}
The graph of $\delta_2(\mathcal{C}^r(x,7))$ are described in Fig.~9.a. Analyzing the above density function we get that the
maximal density is achieved at the endpoint of the above interval (see Fig.~9.b) 
with density $\delta_2(\mathcal{C}^r(w(7)-t(7),7)) \approx 0.81542$. The density in the starting point of the above interval is
$\delta_2(\mathcal{C}^r(0,7)) \approx 0.64805$.
\begin{equation}
\begin{gathered}
\delta_3(\mathcal{C}^r(x,7))=\frac{4 Vol(\mathcal{H}^{h(7)+x} \cap \mathcal{C}^r(7))+4 \cdot Vol(\mathcal{H}^{h(7)-x}
\cap \mathcal{C}^r(7))}{Vol({\mathcal{C}^r(7)})}, \\
\text{where} ~ x \in [0,\min\{h(7),s(7)-h(7)\} \approx 0.41108].
\end{gathered} \tag{3.20}
\end{equation}
\begin{figure}[ht]
\centering
\includegraphics[width=12cm]{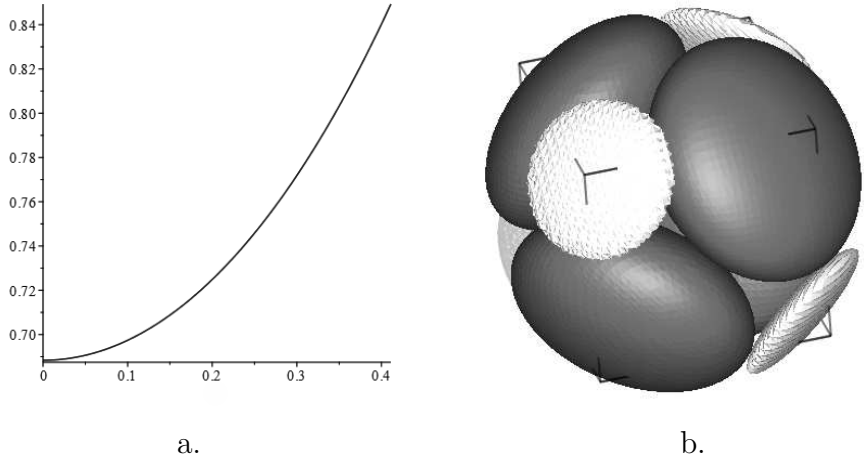}

a. \hspace{6cm} b.
\caption{a.~Density function $\delta_3(x,7))$ $(x\in [0, s(7)-h(7)],7)$. If $s(7)-h(7) \approx 0.41108$ then the density is $\approx 0.84931$. 
The density in the starting point of the above interval is $\delta_3(\mathcal{C}^r(0,7)) \approx 0.68839$.~b.~The non-congruent hyperball arrangement with parameters
$x=s(7)-h(7) \approx 0.41108$, $p=7$ with density $\approx 0.84931$.
}
\label{}
\end{figure}
The graph of $\delta_3(\mathcal{C}^r(x,7))$ are described in Fig.~10.a. Analyzing the above density function we get that the
maximal density is achieved at the endpoint of the above interval (see Fig.~10.b) 
with density $\delta_3(\mathcal{C}^r(s(7)-h(7)\approx 0.41108,7)) \approx 0.84931$. The density in the starting point of the above interval is
$\delta_3(\mathcal{C}^r(0,7) \approx 0.68839$.
\begin{figure}[ht]
\centering
\includegraphics[width=12cm]{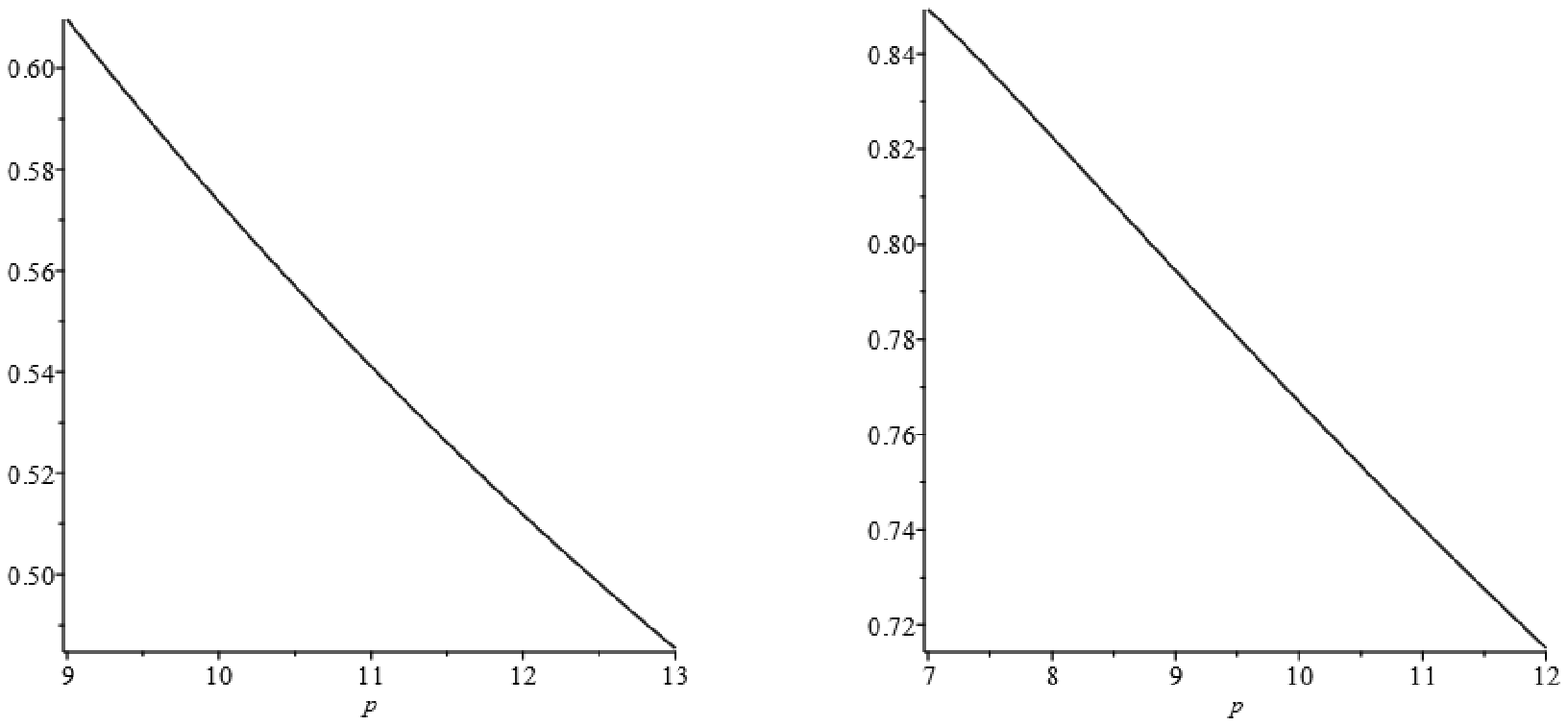}

a. \hspace{6cm} b.
\caption{a.~Density function $\delta_1(\mathcal{C}^r(0,p))$ where $p \in [9, 13]$.
b.~Density function $\delta_3(\mathcal{C}^r(s(7)-h(7) \approx 0.41108,7))$ where $p \in [7,12]$.}
\label{}
\end{figure}
\item Similarly to the above discussions we obtain that if $p=8$ than the maximal density belongs to $\delta_3(\mathcal{C}^r(x,8))$ ($x \in [0, s(p)-h(p) \approx 0.45994]$. 
Analyzing the above density function we get that the maximal density is achieved at the endpoint of the above interval  
with density $\delta_3(\mathcal{C}^r(s(7)-h(p)\approx 0.45994,8)) \approx 0.82259$. 
\item If $p > 8$ then $2h(p) < t(p) < w(p)$ therefore in this case we have to only examine the density functions $\delta_i(\mathcal{C}^r(x,p))$ ($i \in \{1,3\}$). 
Similarly to the above computations for parameter $p=7$ we can analyze the density functions and their maximums of
non-congruent hyperball packings generated by considered truncated cube tilings (or Coxeter tilings $\{4,3,p\}$) for all possible integer parameters
$p > 8$. Using the results of Theorem 3.8 and Corollary 3.9 we obtain the following
\end{enumerate}
\begin{theorem}
\begin{enumerate}
\item The maximum of the density function $\delta_1(\mathcal{C}^r(x,p))$ are attained at the starting
point of the corresponding interval $x \in [0,$ $\min\{h(p),t(p)-h(p)]$ depending on the given integer parameter $p\ge 7$ i.e.
the congruent hyperball packing provides the densest hyperball packing for a given parameter $p$.

(Fig.~11.a shows the graph of the strictly decreasing function $\delta_1(\mathcal{C}^r(0,p))$ if $p\in [9,13]$).
\item If $p=7,8$ then the maximum of the density function $\delta_2(\mathcal{C}^r(x,p))$ is achieved at the endpoint of the interval $[0, 2h(p)-t(p)]$.
If $p>8$ then this case does not occur because $2h(p) < t(p) < w(p)$.
\item The maximum of the density function $\delta_3(\mathcal{C}^r(x,p))$ are attained at the endpoint
point of the corresponding intervals $x \in [0, s(p)-h(p)]$ where $p\ge 7$ is a given integer parameter.

(Fig.~11.b shows the graph of the strictly decreasing function $\delta_3(\mathcal{C}^r(s(p)-h(p),p))$ if $p\in [7,12]$).
\end{enumerate}
\end{theorem}
\begin{theorem}
The maximum of the density functions $\delta_i(\mathcal{C}^r(x,p))$ ($p \ge 7$, integer parameter) ($i \in \{1,2,3\}$) is achieved at the parameters $x=s(p)-h(p)\approx 0.41108$, $p=7$.
Therefore, the density upper bound of the congruent and non-congruent hyperball packings related to the truncated cube tilings $\{4,3,p\}$ $( \mathbb{N} \ni p\ge 7$ is $\approx 0.84931$.
\end{theorem}

\subsubsection{On non-extendable non-congruent hyperball packings $(6 < p < 7, ~ p\in \bR)$}
The computation method described in the former sections is suitable to determine the densities of congruent and non-congruent hyperball packings related 
to the hyperbolic cubes with parameters $(6 < p < 7, ~ p\in \mathbb{R})$.
To any parameter belongs to the truncated cubes we can determine similarly to the
above cases the corresponding densities of their optimal hyperball packings. But these packings cannot be extended to the 3-dimensional
space. Analysing these non-extendable packings for parameters $(6 < p < 7, ~ p\in \mathbb{R})$ we obtain the following
\begin{theorem}
\begin{enumerate}
\item The maximum of the density function $\delta_1(\mathcal{C}^r(x,p))$ $(6 < p < 7, ~ p\in \mathbb{R})$ is attained at the starting
point of the corresponding interval $x \in [0,t(p)-h(p)]$. i.e.
the congruent hyperball packing provides the densest hyperball packing for a given parameter $p$.
(Fig.~12.a shows the graph of the function $\delta_1(\mathcal{C}^r(0,p))$ if $p\in (6,7)$.) This function is attained its maximum at $p^1_{opt} \approx 6.33962$ where
$\delta_1(\mathcal{C}^r(0,p^1_{opt}))\approx 0.70427$.
\item The maximum of the density function $\delta_2(\mathcal{C}^r(x,p))$ $(6 < p < 7, ~ p\in \mathbb{R})$ is achieved at the endpoint of the interval $[0, w(p)-t(p)]$.

(Fig.~12.b shows the graph of the function $\delta_2(\mathcal{C}^r(w(p)-t(p),p))$ if $p\in (6,7)$.) This function is attained its maximum at $p^2_{opt} \approx 6.10563$ where
$\delta_2(\mathcal{C}^r(0,p^2_{opt}))\approx 0.85684$.
\item The maximum of the density function $\delta_3(\mathcal{C}^r(x,p))$ is attained at the endpoint
point of the corresponding interval $x \in [0, s(p)-h(p)]$ where $6 <p < 7$ is a given real parameter.
(Fig.~12.c shows the graph of the function $\delta_3(\mathcal{C}^r(s(p)-h(p),p))$ if $p\in (6,7)$.) This function is attained its maximum at $p^3_{opt} \approx 6.26384$ where
$\delta_3(\mathcal{C}^r(s(p^3_{opt})-h(p^3_{opt}),p^3_{opt}))\approx 0.86145$.
\end{enumerate}
\end{theorem}
\begin{theorem}
The maximum of the density functions $\delta_i(\mathcal{C}^r(x,p))$ ($(6 < p < 7, ~ p\in \mathbb{R})$) ($i \in \{1,2,3\}$) is achieved at the parameters $x=s(p)-h(p)\approx 0.36563$, 
$p_{opt}\approx 6.26384$.
Therefore, the density upper bound of the congruent and non-congruent hyperball packings related to the truncated cube $\{4,3,p\}$  $(6 < p < 7, ~ p\in \mathbb{R})$ 
is $\approx 0.86145$. 
\end{theorem}
\begin{figure}[ht]
\centering
\includegraphics[width=13cm]{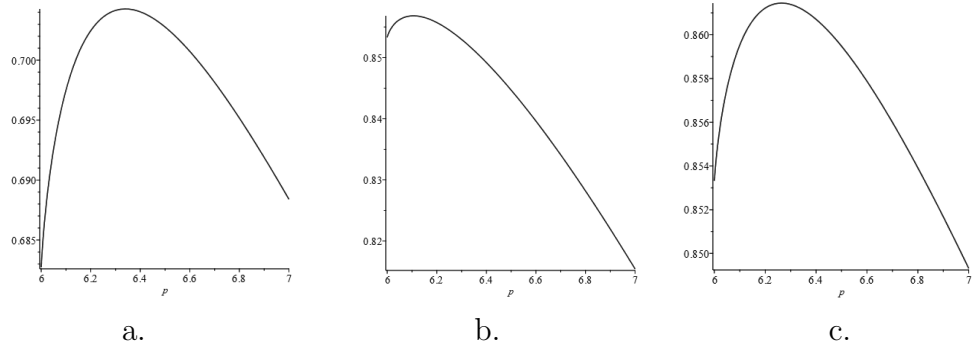}

a. \hspace{4cm} b. \hspace{4cm} c.
\caption{a.~The graph of the function $\delta_1(\mathcal{C}^r(0,p))$ if $p\in (6,7)$. This function is attained its maximum at $p^1_{opt} \approx 6.33962$ where
$\delta_1(\mathcal{C}^r(0,p^1_{opt}))\approx 0.70427$.
b.~ The graph of the function $\delta_2(\mathcal{C}^r(w(p)-t(p),p))$ if $p\in (6,7)$. This function is attained its maximum at $p^2_{opt} \approx 6.10563$ where
$\delta_2(\mathcal{C}^r(0,p^2_{opt}))\approx 0.85684$.
c.~ The graph of the function $\delta_3(\mathcal{C}^r(s(p)-h(p),p))$ if $p\in (6,7)$. This function is attained its maximum at $p^3_{opt} \approx 6.26384$ where
$\delta_3(\mathcal{C}^r(s(p)-h(p),p^3_{opt}))\approx 0.86145$.}
\label{}
\end{figure}
\begin{rmrk}
\begin{enumerate}
\item In our case $\lim_{p\rightarrow 6}(\delta_i(\mathcal{C}^r(x,p)))$ $(i\in\{2,3\}$) is equal to the B\"oröczky-Florian
upper bound of the ball and horoball packings in $\HYP$ (see \cite{B--F64}).
\item The locally optimal hyperball configurations $\delta_2(\mathcal{C}^r(0,p^2_{opt} \approx 6.10563)) 
\approx 0.85684$ and $\delta_3(\mathcal{C}^r(s(p^3_{opt})-h(p^3_{opt}),p^3_{opt}\approx 6.26384))\approx 0.86145$ provide larger densities that the B\"or\"oczky-Florian density 
upper bound $\delta_{BF} \approx 0.85328$ for ball and horoball packings (\cite{B--F64}) but these hyperball packing configurations
are only locally optimal and cannot be extended to the entire hyperbolic 
space $\mathbb{H}^3$.
\end{enumerate}
\end{rmrk}
The problem of finding the densest hyperball (hypersphere) packing with congruent or non-congruent hyperballs
in $n$-dimensional hyperbolic space ($n\ge3$) is not settled yet. 
At this time in $3$-dimensional hyperbolic space $\HYP$ the densest hyperball packing with congruent hyperballs is derived by the regular truncated 
tetrahedron tiling with density $\approx 0.82251$ and with non-congruent hyperballs is derived by the regular truncated 
cube tiling $\{4,3,7\}$ with density $\approx 0.84931$.

But, as we have seen, locally there are hyperball packings with larger density than the
B\"or\"oczky-Florian density upper bound for ball and horoball packings.

We note here, that the discussion of the densest packings in the $n$-dimensional hyperbolic space $n \ge 3$ with horoballs
of different types has not been settled yet as well (see e.g. \cite{KSz}, \cite{KSz14}, \cite{N16}, \cite{Sz12}, \cite{Sz17}). 

\noindent
\footnotesize{Budapest University of Technology and Economics Institute of Mathematics, \\
Department of Geometry, \\
H-1521 Budapest, Hungary. \\
E-mail:~szirmai@math.bme.hu \\
http://www.math.bme.hu/ $^\sim$szirmai}

\end{document}